\documentclass{article}
\usepackage{amsmath}
\usepackage{amsfonts}
\usepackage{amssymb}
\usepackage{amsthm}
\usepackage{amscd}
\usepackage{makeidx}
\makeindex
\newtheorem{theorem}{Theorem}
\newtheorem{remark}{Remark}
\newtheorem{proposition}{Proposition}
\newtheorem{definition}{Definition}

\newcommand{\N}{{\mathbb N}}

\newcommand{\Q}{{\mathcal Q}}
\newcommand{\Ll}{{\mathcal L}}
\newcommand{\Dd}{{\mathcal D}}
\newcommand{\I}{{\mathcal I}}
\begin{document}

\title{On bounded and unbounded idempotents whose sum is a multiple of the 
identity}
\author{Yuri\u\i{}  Samo\u\i{}lenko \vspace{0.3cm}\\
{\footnotesize \sl Institute of Mathematics, Tereshchenkivska 3,  }\\ 
{\footnotesize \sl 252601, Kiev, Ukraine}\\
{\footnotesize \sl e-mail: yurii\_sam@imath.kiev.ua}\\ 
\vspace{0.5cm}\\
Lyudmila  Turowska \vspace{0.3cm}\\
{\footnotesize \sl Department of Mathematics, Chalmers University of 
Technology,}\\ 
{\footnotesize \sl SE-412 96 G\"oteborg, Sweden}\\
{\footnotesize \sl e-mail: turowska@math.chalmers.se}
}
\date{}
\maketitle{}
\begin{abstract}
We study bounded and unbounded representations of the $*$-algebra
$\Q_{n,\lambda}(*)$ generated by $n$ idempotents whose sum equals $\lambda e$
($\lambda\in{\mathbb C}$, $e$ is the identity).
\end{abstract}
\section{Introduction}
Let $H$ be a complex (finite or infinite dimensional) separable Hilbert space 
and $B(H)$ the algebra of bounded operators on $H$.
An operator $Q\in B(H)$ is said to be idempotent if $Q^2=Q$. A selfadjoint 
idempotent $P$ is called orthoprojection.
The problem to study families of idempotent is not just an interesting 
algebraic problem but it also arises
in numerous applications  to  analysis (see, for example, \cite{bgkkrss} and
references therein), 
theory of operators (see \cite{Wu} and references therein), mathematical 
physics
(see \cite{bax,ek} and references therein) etc. 

The structure of pairs of idempotents on finite-dimensional space or, equivalently, 
finite-dimensional representations of the algebra 
${\mathcal Q}_2={\mathbb C}\langle q_1,q_2,e\mid q_k^2=q_k,
k=1,2\rangle$ with the unit $e$, is non-trivial but well understood 
(see \cite{naz, gp}).
The structure of three and more idempotents or, equivalently,
representations of the algebras
${\mathcal Q}_n={\mathbb C}\langle q_1,\ldots,q_n,e\mid q_k^2=q_k,k=1,
\ldots,n\rangle$
for $n\geq 3$ is extremely difficult (the corresponding algebra is wild (see
\cite{df})). Families of idempotents  with
additional condition that their sum is a multiple of the identity
operator or, equivalently, representations of
the algebras ${\mathcal Q}_{n,\lambda}=
{\mathbb C}\langle q_1,\ldots,q_n,e\mid\sum_{k=1}^nq_k=\lambda e,
 q_k^2=q_k,k=1,\ldots,n\rangle$, $\lambda\in{\mathbb C}$, are  studied less
(for  some results about finite-dimensional representations and representations
by bounded operators on Hilbert space of $\Q_{n,\lambda}$ see the first 
subsection at each section).

In this paper we consider the problem of describing up to unitary equivalence
families of idempotents  or, which is equivalent,
$*$-representations of  $*$-algebras generated by
idempotents.

Natural $*$-analogues of the algebras ${\mathcal Q}_n$ and 
${\mathcal Q}_{n,\lambda}$, $n\in{\mathbb N}$, $\lambda\in{\mathbb C}$, are:

(a) the $*$-algebras
${\mathcal Q}_n(*)={\mathbb C}\langle q_1,\ldots,q_n, q_1^*,\ldots,q_n^*,e
\mid q_k^2=q_k,
k=1,\ldots,n\rangle$ generated by $n$ free idempotents and their adjoints.
Denoting ${\mathcal Q}_n^*$ the algebra (not a $*$-algebra) 
${\mathbb C}\langle  q_1^*,\ldots,q_n^*\mid (q_k^*)^2=q_k^*,
k=1,\ldots,n\rangle$, we have ${\mathcal Q}_n(*)={\mathcal Q}_n\star 
{\mathcal Q}_n^*$, where $\star$ is
the sign of free product of algebras;

(b) the $*$-algebras 
${\mathcal Q}_{n,\lambda}(*)=
{\mathbb C}\langle q_1,\ldots,q_n,q_1^*,\ldots,q_n^*,e
\mid \sum_{k=1}^nq_k=\lambda e, q_k^2=q_k,
k=1,\ldots,n\rangle$, $\lambda\in{\mathbb C}$, which are generated by $n$ idempotents
whose sum is a multiple of the unit, and elements which are adjoint to them.
Setting 
${\mathcal Q}_{n,\lambda}^*=
{\mathbb C}\langle q_1^*,\ldots,q_n^*,e
\mid \sum_{k=1}^nq_k^*=\lambda e, (q_k^*)^2=q_k^*,
k=1,\ldots,n\rangle$, we have ${\mathcal Q}_{n,\lambda}(*)=
{\mathcal Q}_{n,\lambda}\star{\mathcal Q}_{n,\lambda}^*$;

and their factor-$*$-algebras, such as

(c) the $*$-algebras 
${\mathcal P}_n={\mathbb C}\langle p_1,\ldots,p_n,e\mid p_k^2=p_k=p_k^*,
k=1,\ldots,n\rangle$ generated by $n$ orthoprojections;

(d) the $*$-algebras
${\mathcal P}_{n,\alpha}={\mathbb C}\langle p_1,\ldots,p_n,e\mid
\sum_{k=1}^np_k=\alpha e, p_k^2=p_k=p_k^*,
k=1,\ldots,n\rangle$, $\alpha\in {\mathbb R}$, 
 generated by $n$ orthoprojections whose sum is a multiple of the unit.

The structure of pairs of orthogonal projections $P_1$, $P_2$, or
representations of $*$-algebra ${\mathcal P}_2$ is well known 
(see, for example,\cite{ras}):
irreducible representations of ${\mathcal P}_2$ are one or two-dimensional
and any representation is a direct sum (or direct integral) of 
irreducible ones.

The structure of three and more orthoprojections or representations
of the $*$-algebras ${\mathcal P}_n$, $n\geq 3$, is very difficult (the $*$-algebras
are $*$-wild). For the definition of $*$-wild algebras we refer
the reader to \cite{krusam,os}.

A number of articles (see \cite{rs,krs}) are devoted to the structure of families of 
orthogonal projections  whose sum is a multiple of the identity or,
equivalently, representations of the $*$-algebras ${\mathcal P}_{n,\alpha}$,
$\alpha\in{\mathbb R}$. In particular,  there were  described
the sets  of $\alpha\in{\mathbb R}$ such that
there exist orthogonal projections $P_1,\ldots P_n$ on a Hilbert space
$H$ so that $\sum_{k=1}^nP_k=\alpha I$, $I$ is the identity operator.
Note that orthoprojections are necessarily bounded operators and $*$-algebras
${\mathcal P}_n$ and ${\mathcal P}_{n,\alpha}$ do not have representations
by unbounded operators.

 Families of idempotents $Q_1,\ldots,Q_n,Q_1^*,\ldots,Q_n^*$ or, equivalently,
representations of $*$-algebra ${\mathcal Q}_n(*)$ have a simple structure 
in the case $n=1$. The situation is similar to the situation for
pairs of orthoprojections (see section 2.2).
If $n\geq 2$ the problem of unitary classification of all families
of idempotents becomes difficult (the $*$-algebra ${\mathcal Q}_{n}(*)$ is
$*$-wild) (\cite{krusam,os}).

In the present paper we study, up to unitary equivalence, the structure of
idempotents whose sum is a multiple of identity, or, equivalently,
representations of the $*$-algebras ${\mathcal Q}_{n,\lambda}(*)$, $\lambda\in
{\mathbb C}$. In contrast to orthogonal projections there exist 
unbounded  operators $Q$ satisfying $Q^2=Q$ (unbounded idempotents) 
(see \cite{popovich}). We study representations
of ${\mathcal Q}_{n,\lambda}(*)$ by bounded operators (i.e., $*$-homomorphisms 
of ${\mathcal Q}_{n,\lambda}(*)$ to $B(H)$) and the sets
$\Lambda_{n,bd}=\{\lambda\in{\mathbb C}\mid \exists Q_1,\ldots,Q_n\in B(H),
\sum_{k=1}^nQ_{k}=\lambda I, Q_k^2=Q_k, k=1,\ldots,n\}$ together with
representations
by unbounded operators and the corresponding sets $\Lambda_{n,unbd}$
(for  exact definitions see subsection 3 at each section).
 
For $n=3$ (sect.~3.2) representations of ${\mathcal Q}_{n,\lambda}(*)$ by 
bounded operators  exist only for $\lambda\in\{0,1,3/2,2,3\}$, the $*$-algebras
${\mathcal Q}_{3,0}(*)\simeq{\mathcal Q}_{3,3}(*)$ are one-dimensional, 
the $*$-algebras ${\mathcal Q}_{3,1}(*)\simeq {\mathcal Q}_{3,2}(*)\simeq
{\mathbb C}^3\star
{\mathbb C}^3$ and 
${\mathcal Q}_{3,3/2}(*)\simeq M_2({\mathbb C})\star M_2({\mathbb C})$ are $*$-wild.
  For $*$-wild algebras we also study
additional conditions under which the problem of unitary classification
of their representations become transparent. Properties of the 
representations of ${\mathcal Q}_{n,\lambda}(*)$ by unbounded 
operators are essentially 
the same as representations by bounded operators  in the case $n=3$ 
(sect.\  3.3). 
In particular, $\Lambda_{3,bd}=\Lambda_{3,unbd}$.

For $n=4$ the situation becomes different. For example, the sum of four 
bounded idempotents equals zero only if the idempotents are zero operators, 
which is not the case for unbounded ones (see \cite{bes}). If the sum of
bounded idempotents is equal to $1$ then the idempotents are mutually
orthogonal, but there exist unbounded non-orthogonal idempotents with the sum
equal to $1$ (see \cite{erss}). Moreover, 
$\Lambda_{4,bd}=\{0,1,1+\frac{k}{k+2} 
(k\in{\mathbb N}),2, 3-\frac{k}{k+2},3,4\}$, while
$\Lambda_{4,unbd}={\mathbb C}$ (Proposition~\ref{lambda4unbd}).
In Section~4 we study the problem of describing of representations
of ${\mathcal Q}_{4,\lambda}(*)$ by bounded and unbounded operators and 
 representations
of  ${\mathcal Q}_{4,\lambda}(*)$ under some additional restrictions.
The case $\Q_{4,0}(*)$ is treated in details.
 
In Section~5, following \cite{rs},  we show that $\Lambda_{5,bd}={\mathbb C}$ 
and that the problem of unitary classification of already
bounded representations of 
the $*$-algebra ${\mathcal Q}_{n,\lambda}(*)$ for $n\geq 5$ is difficult: 
the $*$-algebra is not of type $I$ for any $\lambda\in{\mathbb C}$ and 
it is $*$-wild for some $\lambda\in{\mathbb R}$.

We would like to emphasise one more time that speaking about representation
of algebras we mean homomorphisms into algebras of linear operators on a 
vector space and  
$*$-homomorphisms into a $*$-algebra of linear bounded or unbounded operators
defined on a Hilbert space
if we speak about representations or $*$-representations of 
$*$-algebras. We will also restrict 
ourself to indecomposable finite-dimensional representations up to similarity 
when talk about description of representations of algebras.
Description of $*$-representations is reducing to the description of 
irreducible (=indecomposable) representations up to unitary equivalence.


\section{Representations of algebras $\Q_{2,\lambda}$ and $*$-algebras
 ${\mathcal Q}_{2,\lambda}(*)$}

\subsection{Algebras ${\mathcal Q}_{2,\lambda}$ and their representations}

Algebra $Q_{2,\lambda}$ is non-zero only for $\lambda\in\Lambda_2=\{0,1,2\}$.
We have ${\mathcal Q}_{2,0}={\mathcal Q}_{2,2}={\mathbb C}e$.
The algebra ${\mathcal Q}_{2,1}$ is easily seen to be equal to 
${\Q}_1={\mathbb C}\langle q,e|q^2=q\rangle$, its finite-dimensional 
indecomposable representations are  one-dimensional: $\pi(q)=0$ or 
$\pi(q)=1$.

\subsection{$*$-Algebras ${\mathcal Q}_{2,\lambda}(*)$ and their 
representations by bounded operators}

$*$-Representations of $\Q_{2,0}(*)=\Q_{2,2}(*)={\mathbb C}e$ are trivial. 
The problem of unitary classification of  $\Q_{2,1}(*)=\Q_1(*)$ reduces to 
the problem
of describing a single idempotent up to unitary equivalence. Any irreducible 
representation
of $\Q_{1}(*)$ is one- or two-dimensional
and given by

\begin{eqnarray}
&&(a)\quad \pi(q_1)=1,\ \quad \text{ or }\quad \pi(q_1)=0; \label{rep1}\\
&&(b)\quad \pi(q_1)=\left(\begin{array}{cc}
1&y\\
0&0
\end{array}\right), \  y\in(0,\infty).\label{rep2}
\end{eqnarray}
(see, for example \cite{yugoslav,os}).
The structure of arbitrary
representation of $\Q_1(*)$ by bounded operators is given by the following
statement (see \cite{krusam1}):
for any representation $\pi$ on a Hilbert space $H$ there exist a unique
decomposition $H=H_0\oplus H_1\oplus{\mathbb C}^2\otimes H_2$ and 
a projection-valued measure $dE(\cdot)$ on $H_2$ whose support is a bounded
subset of $(0,\infty)$ and such that
\begin{equation}\label{structure_th}
\pi(q_1)=0\cdot I_{H_0}\oplus I_{H_1}\oplus\int_0^{\infty}
\left(\begin{array}{cc}
1&y\\
0&0
\end{array}\right)\otimes dE(y).
\end{equation}
 
\subsection{Representations of ${\mathcal Q}_{2,\lambda}(*)$ by
unbounded operators}
Here we describe $*$-representations of $*$-algebra $\Q_{2,1}(*)=Q_1(*)$ by
unbounded operator recalling necessary definitions of the concepts involved.

Let $\Phi$ be a dense linear subset of a Hilbert space, $H$. Let
${\Ll}^+(\Phi)=\{X\in{\Ll}(\Phi)\mid \Phi\subset \Dd(X^*), 
X^*\Phi\subset\Phi\}$, here $\Dd(X)$ is the domain of the operator $X$.
Then ${\Ll}^+(\Phi)$ is an algebra with involution $X^+=X^*|_{\Phi}$.
By {\it $*$-representation} of a $*$-algebra
${\mathfrak A}$ by unbounded operators we call a unital $*$-homomorphism
$\pi:{\mathfrak A}\to{\Ll}^+(\Phi)$ (see, for example, \cite{inoue,Shbook}).
We write also $D(\pi)$
 for the domain $\Phi$ and call it the domain of the representation $\pi$.

Define $\Lambda_{n,unbd}$ to be the set of all $\lambda\in{\mathbb C}$ such that
there exists a $*$-representation of ${\mathcal Q}_{n,\lambda}(*)$ by
unbounded operators.
Since  $*$-algebras $\Q_{2,\lambda}(*)$ are non-zero only for 
$\lambda\in\Lambda_2$, we have $\Lambda_{2,bd}=\Lambda_{2,unbd}=\Lambda_2$.

The class of representations defined above is very large and practically 
indescribable (see \cite{samtur,Shbook,t}). So 
if one wishes to get structure theorems giving a description of 
$*$-representations up to unitary equivalence then one should impose some
additional conditions on the domain $\Phi$.  For example, one can require 
that $\Phi$ consists of bounded (entire, analytical vectors) for some
operators of the representation. 
Recall that  a vector $f\in \cap_{k\in{\N}}\Dd (X^k)\subset H$ is called {\it bounded}
({\it entire}, {\it analytical}) vector for operator $X$ on $H$ if
there is a constant $c_f$  such that 
$||X^nf||\leq c_f^n||f||$ for any $n\in{\N}$ (the function
$\sum_{k=0}^{\infty}(||X^kf||/k!)z^k$ is entire or analytical at $0$). 
The set of bounded
(entire, analytical) vectors for operator $X$ will be denoted by
$H_b(X)$ ($H_c(X)$ and $H_a(X)$ respectively). 
Imposing this type of conditions
we will call these representations {\it integrable} or {\it well-behaved} or 
{\it good} following
the terminology in the theory of representations of Lie algebras 
(\cite{nel,Shbook}).
Definitions of equivalent representations, irreducible 
representations which are necessary for formulating structure theorems, 
will be given for every particular class of representation considered in the 
paper.

Let $\pi$ be a representation of $\Q_{1}(*)$ defined on a domain $D(\pi)$
of a Hilbert space $H(\pi)$.
We say that $\pi$ is a {\it well-behaved} representation if $\Delta=
\overline{\pi(qq^*+qq^*-(q+q^*))}$
(the closure of the operator $\pi(qq^*+qq^*-(q+q^*))$) is selfadjoint and 
$D(\pi)=H_b(\Delta)$. Note that
$qq^*+qq^*-(q+q^*)$ is a central element  of $\Q_{1}(*)$.
Setting $a=q+q^*-e$ and $b=i(q^*-q)$, which are clearly selfadjoint,
we obtain $qq^*+qq^*-(q+q^*)=a^2+b^2-e$.

Let $H_m(\pi)$ denote the set of all vectors $f\in H_b(\Delta)$ such that
$||\Delta^nf||\leq m^n||f||$. Since $qq^*+qq^*-(q+q^*)$ is central and
$\Delta$ is selfadjoint, one can show that the subspaces $H_m(\pi)$ is 
reducing $\pi$, i.e., $\pi=\pi_1\oplus\pi_2$ and $H_m(\pi)=D(\pi_1)$.
Moreover, each subrepresentation $\pi|_{H_m(\pi)}$ is bounded.
In order to see this, it is enough to show that $\pi(a)$ and $\pi(b)$ are 
bounded on $H_m(\pi)$. Given $f\in H_m(\pi)$, we have
$||\pi(a^2+b^2)f||\leq (m+1)||f||$ and $||\pi(a)f||^2=
(\pi(a)^2f,f)\leq (\pi(a^2+b^2)f,f)\leq ||\pi(a^2+b^2)f||\cdot||f||\leq 
(m+1)||f||^2$, the same holds for $\pi(b)$.
 
A representation $\pi$ is called {\it irreducible} if the only 
linear 
subspace reducing  $\pi$ are $\{0\}$ and $D(\pi)$. Since for a well-behaved
representation $\pi$ at least one of 
$H_m(\pi)$ is non-zero we 
obtain that any  well-behaved irreducible representation is bounded.

We say that two well-behaved representations $\pi_1$ and $\pi_2$ are 
{\it unitarily equivalent}
if there exist a unitary operator $U:H(\pi_1)\to H(\pi_2)$ such that $UH_m(\pi_1)=
H_m(\pi_2)$ and $U\pi_1(a)f=\pi_2(a)Uf$ for any $f\in H_m(\pi_1)$ and each 
$m\in{\mathbb N}$. 

Now we state a structure theorem. Its proof essentially
follows the proof of an analogous statement for bounded representations of
$\Q_1(*)$.

\begin{proposition}\label{idempotent}
Any irreducible well-behaved representation of $\Q_{1}(*)$ is bounded and, up
to unitary equivalence, is  
given by $(\ref{rep1})-(\ref{rep2})$.
Any representation $\pi$ of $\Q_{1}(*)$ is a direct sum (direct integral) of
irreducible ones and given by (\ref{structure_th}), where the equality holds on
vectors $f\in D(\pi)$.
\end{proposition}

\begin{remark}\label{rep_1}\rm
If a $*$-algebra is defined in terms of generators and commutation relations 
(the algebraic equalities imposed on the generators), instead of
$*$-representations of the $*$-algebra one  speaks often about representations of 
this set of relations, i.e. families of operators satisfying the relations
on some invariant dense domain.

Let $Q$, $Q^*$ be closed operators such that that there exists  a dense
domain $\Phi$ satisfying the following conditions: (1) $\Phi$ is invariant 
with respect
to $Q$, $Q^*$; (2) $\Phi$ is a core for $Q$, $Q^*$ (i.e. the closure of 
operators $Q|_{\Phi}$ and $Q^*|_{\Phi}$ are $Q$ and $Q^*$ respectively); (3) 
$\Phi\subset
H_a(QQ^*+Q^*Q-(Q+Q^*))$; (4) $Q^2f=Qf$ for any $f\in\Phi$.
 Then $Q$, $Q^*$ generate a $*$-representation $\pi$ of $\Q_{1}(*)$ with 
the domain $\Phi$ by setting $\pi(q)f=Qf$, $\pi(q^*)f=Q^*f$, $f\in \Phi$  
and then extending it to the whole algebra. 
One can show that the representation $\pi^*$ ($\pi^*(a)=\pi(a^*)^*|_{D(\pi^*)}$,
and $D(\pi^*)=\cap_{a\in Q_{2,1}(*)}
D(\pi(a)^*)$) is a unique selfadjoint $*$-representation 
$\rho$ of $\Q_{1}(*)$ such that $\overline{\rho(q)}=Q$, 
$\overline{\rho(q^*)}=Q^*$. Recall that $\rho$ is selfadjoint if $\rho=\rho^*$.
Note that the domain $\Phi$ satisfying
$(1)-(4)$ is not uniquely defined. But for every choice of $\Phi$ we have
the same selfadjoint representation.  
In particular, the domain $\Phi$ can be chosen to be equal 
$H_b(QQ^*+Q^*Q-(Q+Q^*))$ which implies that
the unique selfadjoint representation $\rho$ is $\pi^*$, where $\pi$ is a 
well-behaved representation. In this case the representation 
$\pi^*$ is irreducible iff
$\pi$ is irreducible, two well-behaved representations $\pi_1$ and $\pi_2$ 
are unitarily equivalent iff
$\pi_1^*$ and $\pi_2^*$ are unitarily equivalent in the sense that there exist
a unitary operator $U$ of $H(\pi_1^*)$ onto $H(\pi_2^*)$ such that
$UD(\pi_1^*)=D(\pi_2^*)$ and $U^{-1}\pi_2^*(a)Uf=\pi_1^*(a)f$, 
$f\in D(\pi_1^*)$. So the problem to classify all well-behaved representation
is equivalent to the problem to classify all selfadjoint representations
defined above or all pairs of operators $(Q,Q^*)$.  
For  concepts of the theory of representations by unbounded 
operators we refer the reader to \cite{Shbook}. 
\end{remark}

\begin{remark}\label{c*alg}\rm
There is also a correspondence between well-behaved rep\-re\-sen\-ta\-tions 
and 
repre\-sen\-ta\-tions arising from representations of some $C^*$-algebra.
Let ${\mathcal A}=\{f\in C([0,\infty), M_2({\mathbb C}))\mid
f(0)\text{ is diagonal}, \lim_{t\to\infty} f(t)=0\}$. ${\mathcal A}$ is a
$C^*$-algebra. Let $q^*(t)=\left(\begin{array}{cc}
1&0\\t&0\end{array}\right)$. One can show that 
$q^*\in C([0,\infty), M_2({\mathbb C}))$
is affiliated with ${\mathcal A}$. Moreover, ${\mathcal A}$ is generated
by $q^*$ and there exists a dense domain $D$ of ${\mathcal A}$ (for example,
$D=\{f\in{\mathcal A}\mid \text{supp }f\text{ is compact}\}$) which is
invariant with respect to $q=(q^*)^*$, $q^*$, $D$ is a core for $q$ and 
$q^*$ and such that $q^2a=qa$, $(q^*)^2a=q^*a$ for any $a\in D$.
For the notion of affiliated elements and $C^*$-algebras generated 
by unbounded elements we refer the reader to \cite{wor1,wor2} (see also section
4.3).
Let ${\mathcal R}$ denote the set of pairs $(Q,Q^*)$  satisfying
the conditions given in Remark~\ref{rep_1}. Then
${\mathcal R}=\{(\pi(q),\pi(q^*))\mid \pi\text { is a non-degenerate
representation of }{\mathcal A}\}$ (recall that $\pi$ is non-degenerate if 
$\overline{\pi({\mathcal A})H}=H$). 
Here $\pi(q)$, $\pi(q^*)$ is the unique extension of the representation $\pi$
to affiliated elements.
This was essentially proved in \cite{popovich}. It means that any
well-behaved representation arises from a representation of a $C^*$-algebra
and any representation of ${\mathcal A}$ gives rise to a well-behaved 
representation of $\Q_1(*)$.
\end{remark} 

\section{Representations of algebras ${\mathcal Q}_{3,\lambda}$ and 
$*$-algebras
${\mathcal Q}_{3,\lambda}(*)$}

\subsection{Algebras ${\mathcal Q}_{3,\lambda}$ and their representations}
All algebras $\Q_{3,\lambda}$ are finite-dimensional. They are non-zero only
if $\lambda\in\Lambda_3=\{0,1,3/2,2,3\}$.
We have $\Q_{3,0}=\Q_{3,3}={\mathbb C}e$ with trivial representations.
The idempotents $q_i$, $i=1,2,3$, of $\Q_{3,1}=\Q_{3,2}$ are orthogonal, i.e.
$q_iq_j=0$ ($i\ne j$). In fact, since $(e-q_3)^2=e-q_3$, we have 
$\{q_1,q_2\}=q_1q_2+q_2q_1=0$. But the idempotents anti-commute iff
$q_1q_2=q_2q_1=0$.  Therefore   
$\Q_{3,1}=\Q_{3,2}={\mathbb C}q_1\oplus{\mathbb C}q_2\oplus {\mathbb C}q_3
= {\mathbb C}^3=\Q_{2,\perp}$ 
($={\mathbb C}\langle q_1,q_2,e\mid q_i^2=q_i, q_1q_2=
q_2q_1=0\rangle$).
The algebra $\Q_{3,3/2}$ is isomorphic to $ M_2({\mathbb C})$. 
Finite-dimensional representations of ${\mathbb C}^3$ and $M_2({\mathbb C})$
are easy to describe.

\subsection{$*$-Algebras ${\mathcal Q}_{3,\lambda}(*)$ and their 
representations by bounded operators}
As in the algebraic situation we have $\Lambda_{3,bd}=\{0,1,3/2,2,3\}$.
In contrast to  $\Q_{3,0}(*)=\Q_{3,3}(*)={\mathbb C}e$ whose representations
are trivial,
the structure of representations of $\Q_{3,1}(*)=\Q_{3,2}(*)$ and 
$\Q_{3,3/2}(*)$ is very 
complicated (the $*$-algebras are $*$-wild).
Before proving this we recall some results and constructions concerning
wild $*$-algebras. For the general definition and results we refer the reader
to \cite{krusam,os}.

Let ${\mathfrak S}_2$ denote the unital $*$-algebra generated by free
selfadjoint elements, $a$, $b$, and let $C^*({\mathcal F}_2)$ be the group
$C^*$-algebra of the free group ${\mathcal F}_2$ with two generators.
For a unital $*$-algebra $A$ we denote by $\text{Rep }A$ the category
of $*$-representations of $A$ whose objects are unital $*$-representations
of $A$ considered up to unitary equivalence and its morphisms are intertwining
operators. If $A$, $B$ are $*$-algebras, $A\otimes B$ denote the $*$-algebra
that consists of all finite sums of the form $a\otimes b$, $a\in A$, $b\in B$.

Assume that for a unital $*$-algebra, ${\mathfrak A}$, there exists
a unital $*$-homomorphism $\psi:{\mathfrak A}\to 
M_n({\mathbb C})\otimes{\mathfrak S}_2$ (or to $M_n({\mathbb C})\otimes
C^*({\mathcal F}_2)$) for some $n\in{\mathbb N}$.
Then $\psi$ generates a functor $F_{\psi}:\text{Rep }{\mathfrak S}_2
\to \text{Rep }{\mathfrak A}$ (or  $F_{\psi}:\text{Rep } C^*({\mathcal F}_2)
\to \text{Rep }{\mathfrak A}$) defined as follows:
\begin{itemize}
\item
$F_{\psi}(\pi)=id\otimes\pi$ for any $\pi\in \text{Rep }{\mathfrak S}_2$ (or
$\pi\in \text{Rep } C^*({\mathcal F}_2)$), where $id$ is the identity 
representation
of $M_n({\mathbb C})$;

\item
$F_{\psi}(C)=I_n\otimes C$ if $C$  intertwines
representations $\pi_1$, $\pi_2\in \text{Rep }{\mathfrak S}_2$ 
(or  $\pi_1$, $\pi_2\in\text{Rep } C^*({\mathcal F}_2)$), here $I_n$ is the 
identity operator on ${\mathbb C}^n$.
\end{itemize}
If the functor $F_{\psi}$ is full then the $*$-algebra ${\mathfrak A}$ is
$*$-wild. To see that $F_{\psi}$ is full one has to check that  any operator 
intertwining two
representations $F_{\psi}(\pi_1)$, $F_{\psi}(\pi_2)$ with 
$\pi_1$, $\pi_2\in \text{Rep }{\mathfrak S}_2$ (or  $\pi_1$, $\pi_2\in
\text{Rep } C^*({\mathcal F}_2)$) is equal to
$F_{\psi}(C)$, where $C$ is an operator which intertwines $\pi_1$ 
and $\pi_2$.
In particular, we have  that $F_{\psi}(\pi)$, 
 is irreducible iff $\pi$ is irreducible, two
representations $F_{\psi}(\pi_1)$, $F_{\psi}(\pi_2)$ are unitarily
equivalent iff $\pi_1$ and $\pi_2$ are unitarily equivalent.
In this case we say that the problem of unitary classification of
all representations of ${\mathfrak A}$ contains as a subproblem the problem
of unitary classification of representations
of ${\mathfrak S}_2$ (or $C^*({\mathcal F}_2)$).
 
\begin{proposition}\label{wild1}
(a) The $*$-algebra $\Q_{3,1}(*)=\Q_{3,2}(*)=\Q_{2,\perp}(*)$ is $*$-wild.

(b) The $*$-algebra $\Q_{3,3/2}(*)$ is $*$-wild.
\end{proposition}

\begin{proof}
(a) Following \cite[Theorem~6]{krusam} or \cite[Theorem~59]{os}, 
define a $*$-homomorphism
$\psi:\Q_{2,\perp}(*)\to M_3({\mathbb C})\otimes {\mathfrak S}_2$ by
$$
\psi(q_1)=\left(\begin{array}{ccc}
e&e&a+ib\\
0&0&0\\
0&0&0
\end{array}\right), \ \psi(q_2)=\left(\begin{array}{ccc}
0&-e&-e\\
0&e&e\\
0&0&0
\end{array}\right).$$

One can easily check that the generated functor $F_{\psi}$ is full,
i.e., the $*$-algebra $\Q_{2,\perp}(*)=\Q_{3,1}(*)$ is
$*$-wild.

(b) To see that $\Q_{3,3/2}(*)\simeq M_2({\mathbb C})\star
M_2({\mathbb C})^*$ is  $*$-wild we define a $*$-homomorphism
$\psi: 
M_2({\mathbb C})\star M_2({\mathbb C})^*\to M_2({\mathbb C})\otimes 
{\mathfrak S}_2$ by
$$\psi(e_{11})=\left(\begin{array}{cc}
e&-a-ib)\\
0&0
\end{array}\right),\ \psi(e_{12})=\left(\begin{array}{cc}
0&e\\
0&0
\end{array}\right), \psi(e_{21})=\left(\begin{array}{cc}
a+ib&-(a+ib)^2\\
e&-a-ib)
\end{array}\right),$$
where $e_{ij}$ are the matrix units in $M_2({\mathbb C})$.
We leave it to the reader to check that
the functor $F_{\psi}$ is full, i.e., the $*$-algebra
$\Q_{3,3/2}(*)\simeq M_2({\mathbb C})\star
M_2({\mathbb C})^*$ is $*$-wild.
\end{proof}
Note that by Proposition~\ref{wild1} the $*$-algebras $Q_{2,\perp}(*)$ and 
$Q_{3,3/2}(*)$ have infinite-dimensional
irreducible representations.

Imposing additional conditions on representations, the problem of describing
them up to unitary equivalence might become easier.
For example, representations of ${\mathcal P}_{3,1}={\mathcal P}_{3,2}$ and 
${\mathcal P}_{3,3/2}$ or, 
equivalently, representations of $\Q_{3,1}(*)=\Q_{3,2}(*)=\Q_{2,\perp}(*)$ and 
$\Q_{3,3/2}(*)$
with the condition that the images of $q_i$, $i=1,2,3$ are selfadjoint are very
simple. There exist three non-unitarily equivalent irreducible
representations of ${\mathcal P}_{3,1}$: $\pi_i(p_i)=1$, $\pi_i(p_j)=0$,
$i\ne j$, $i=1,2,3$, on $H={\mathbb C}^1$. There exists a unique irreducible 
representation of  ${\mathcal P}_{3,3/2}$ on $H={\mathbb C}^2$:
$$\pi(p_1)=\left(\begin{array}{cc}
1&0\\
0&0\end{array}\right),\ 
\pi(p_2)=\left(\begin{array}{cc}
1/4&\sqrt{3}/4\\
\sqrt{3}/4&3/4\end{array}\right),\
\pi(p_3)=\left(\begin{array}{cc}
1/4&-\sqrt{3}/4\\
-\sqrt{3}/4&3/4\end{array}\right).$$

Requiring that images of $q_1,q_2\in\Q_{2,\perp}(*)$ satisfy the conditions
$\pi(q_1)\pi(q_2^*)=0$ and $\pi(q_2^*)\pi(q_1)=0$ we obtain that  irreducible
representations are one- or two-dimensional and, up to unitary equivalence,
given by

\begin{equation}\label{q_21}
\begin{array}{c}
\text{(a) } \pi(q_1)=\varepsilon_1, \ \pi(q_2)=\varepsilon_2,\ \varepsilon_i\in\{0,1\},
\varepsilon_1\varepsilon_2=0;\\
\text {(b) } \pi(q_1)=\left(\begin{array}{cc}
1&\alpha\\
0&0\end{array}\right), \ \pi(q_2)=0 \text{ and } 
\pi(q_1)=0, \ \pi(q_2)=\left(\begin{array}{cc}
1&\alpha\\
0&0\end{array}\right),
\text{ where } \alpha>0.
\end{array}
\end{equation}

\subsection{Representations of ${\mathcal Q}_{3,\lambda}(*)$ by
unbounded operators}
Since the $*$-algebra $\Q_{3,\lambda}(*)$ is non-zero only for 
$\lambda\in\Lambda_3=\Lambda_{3,bd}$, we have 
$\Lambda_{3,unbd}=\Lambda_3=\{0,1,3/2,2,3\}$. 

Consider unbounded rep\-re\-sen\-ta\-tions of a fac\-tor-$*$-algeb\-ra
of $\Q_{3,1}(*)=\Q_{2,\perp}(*)$, namely, $\Q_{2,\perp}(*)/J$, where $J$ is 
the tw-sided $*$-ideal 
generated by $q_1q_2^*$ and $q_2^*q_1$.
Let $\pi$ be a representation of $\Q_{2,\perp}(*)/J$ on a domain $D(\pi)\subset H(\pi)$ and let 
$$\Delta_1=\overline{\pi(q_1q_1^*+q_1^*q_1-(q_1+q_1^*))},\  
\Delta_2=\overline{\pi(q_2q_2^*+q_2^*q_2-(q_2+q_2^*))}.$$
We  say that $\pi$ is {\it well-behaved} if $\Delta_1$ and $\Delta_2$ are selfadjoint, 
$\Delta_1$, $\Delta_2$ strongly commute (i.e., spectral projections of
$\Delta_1$, $\Delta_2$ mutually commute)
 and $D(\pi)=H_b(\Delta_1,\Delta_2)$, the set of bounded
vectors for both $\Delta_1$ and $\Delta_2$. For each $m\in{\mathbb N}$ denote 
by $H_m(\pi)$ the set 
$E_{\Delta_1}((-m,m))E_{\Delta_2}((-m,m))H(\pi)$, where $E_{\Delta_i}(\cdot)$
is a resolution of the identity for the selfadjoint operator $\Delta_i$.
We have $D(\pi)=\cup_{m\in{\mathbb N}} H_m(\pi)$.
It is easy to see that $H_m(\pi)$ is reducing $\pi$ and 
the subrepresentation 
$\pi|_{H_m(\pi)}$ is bounded. With the same definition of irreducibility and
unitary equivalence as in section~2.3 we have 
\begin{proposition}
Any irreducible well-behaved representation $\pi$ of $\Q_{2,\perp}(*)/J$
is bounded and given by (\ref{q_21}). 

For any well-behaved representation
$\pi$ on a Hilbert space $H(\pi)$  there exist a unique decomposition 
$H(\pi)=H_{00}\oplus H_{01}\oplus H_{10}\oplus {\mathbb C}^2\otimes H_{1}
\oplus {\mathbb C}^2\otimes H_2$ and projection-valued measures
$dE_1(\cdot)$, $dE_2(\cdot)$ on $H_1$ and $H_2$ respectively such that
\begin{eqnarray*}
\pi(q_1)=0\cdot I_{H_{00}}\oplus 0\cdot I_{H_{01}}\oplus I_{ H_{10}}
\oplus\int_0^{\infty}
\left(\begin{array}{cc}
1&y\\
0&0
\end{array}\right)\otimes dE_1(y)\oplus
0\cdot I_{H_2},\\
\pi(q_2)=0\cdot I_{H_{00}}\oplus I_{H_{01}}\oplus 0\cdot I_{H_{10}}
\oplus 0\cdot I_{H_1}\oplus\int_0^{\infty}
\left(\begin{array}{cc}
1&y\\
0&0
\end{array}\right)\otimes dE_2(y),
\end{eqnarray*}
where the equalities hold on $D(\pi)$.
\end{proposition}
\begin{proof}
The first statement follows from the fact that any irreducible representation
is bounded. The second one follows from Proposition~\ref{idempotent}. In fact,
since $$\Delta_2\pi(q_i)f=\pi(q_i)\Delta_2f, \ \Delta_2\pi(q_i^*)f=\pi(q_i^*)\Delta_2f$$
for any $f\in D(\pi)=H_b(\Delta_1,\Delta_2)$, $i=1,2$, it follows 
that
 $$E_{\Delta_2}(\delta)\pi(q_i)f=\pi(q_i)E_{\Delta_2}(\delta)f,\
E_{\Delta_2}(\delta)\pi(q_i^*)f=\pi(q_i^*)E_{\Delta_2}(\delta)f, \ i=1,2,$$
 for any Borel  $\delta$ and any $f\in D(\pi)$ (see \cite{osb}[Theorem~1]).
This implies that $E_{\Delta_2}(\{0\})H(\pi)$ is reducing $\pi$: $H(\pi)=H^1\oplus H^2$,
where $H^1=E_{\Delta_2}(\{0\})H(\pi)$, $H^2=H_1^{\perp}$, and $\pi=\pi_1\oplus\pi_2$,
where $\pi_i=\pi|_{D(\pi)\cap H^i}$, $i=1,2$.

Since $\Delta_2|_{H_1}=0$, we get that $\pi_1(q_2)$, $\pi_1(q_2^*)$ are 
bounded. 
Therefore, restricting $\pi_1$ to the $*$-subalgebra 
${\mathbb C}\langle q_2,q_2^*,e\mid q_2^2=q_2\rangle=
\Q_1(*)$, we obtain a bounded representation of $\Q_1(*)$ such that
the image of $q_2q_2^*+q_2^*q_2-(q_2+q_2^*)$ is zero. It follows from  
structure
theorem for bounded representations of $\Q_1(*)$ (section~2.2) that 
$\pi_1(q_2)$ is an 
orthoprojection.
Using the same arguments as before and the fact that $\pi(q_1)\pi(q_2)=
\pi(q_2)\pi(q_1)=\pi(q_1^*)\pi(q_2)=\pi(q_2)\pi(q_1^*)=0$, one obtains that
$H_1^1=\ker\pi(q_2)$ is reducing $\pi_1$, $\pi_1^1=\pi_1|_{H_1^1}$ restricted
to the $*$-subalgebra ${\mathbb C}\langle q_1,q_1^*,e\mid q_1^2=q_1\rangle=
\Q_1(*)$ is a well-behaved representation of $\Q_1(*)$ in the sense given 
in section
2.3, $\pi_1(q_1)|_{(H_1^1)^{\perp}}=0$, 
$\pi_1(q_2)|_{(H_1^1)^{\perp}}=I$.

Since the  kernel of $\Delta_2|_{H^2}$ is $\{0\}$ and 
$\pi(q_i)\pi(q_j)=
\pi(q_i)\pi(q_j^*)=\pi(q_j^*)\pi(q_i^*)=0$, $i\ne j$, we obtain
$\pi_2(q_1)=\pi_2(q_1^*)=0$ 
 and $\pi_2$ restricted to $*$-subalgebra
${\mathbb C}\langle q_2,q_2^*,e\mid q_2^2=q_2\rangle=
\Q_1(*)$ is a well-behaved representations of $\Q_1(*)$. This completes 
the proof.
\end{proof}
Clearly, if the support of the measures $dE_1(\cdot)$ and $dE_2(\cdot)$ are
unbounded, the representation $\pi$ from the proposition is unbounded.
Any well-behaved representation of $\Q_{2,\perp}(*)/J$ is a representation of
$Q_{2,\perp}(*)$.   
\begin{remark}\rm
As for $\Q_1(*)$ there exists a correspondence between well-behaved
representations of $\Q_{2,\perp}(*)/J$ and representations of some 
$C^*$-algebra, namely the $C^*$-algebra ${\mathcal A}\oplus{\mathcal A}$, where
${\mathcal A}=\{f\in C([0,\infty),M_2({\mathbb C}))\mid f(0)
\text{ is diagonal }, \lim_{t\to\infty} f(t)=0\}$. Setting
$q_1(t)=q(t)\oplus 0$, $q_2(t)=0\oplus q(t)$, where $q\ \eta\ {\mathcal A}$ is
the element which was defined in Remark~\ref{c*alg}, we obtain that
$q_1$, $q_2$ generate  ${\mathcal A}\oplus{\mathcal A}$
as affiliated elements, there exists a dense domain 
$ D'\subset{\mathcal A}\oplus{\mathcal A}$ (for example,
$ D'=D\oplus D$, where $D=\{f\in{\mathcal A}\mid \text{supp }f\text{ is compcat}\}$) such that $D'$ is invariant
with respect to $q_i$, $q_i^*$, $i=1,2$, $D$ is a core for $q_i$, $q_i^*$
and  the relations $q_i^2=q_i$, $q_1q_2=q_2q_1=q_1^*q_2=q_2q_1^*=0$ hold
on $D'$. Moreover, if we let ${\mathcal R}$ denote the set of pairs
$(\overline{\pi(q_1)},\overline{\pi(q_2)})$, where $\pi$ is a well-behaved
 representation of $\Q_{2,\perp}(*)/J$ 
then ${\mathcal R}=\{\rho(q_1),\rho(q_2)
\mid \rho\text{ is a non-degenerate representation of } 
{\mathcal A}\oplus{\mathcal A}\}$. Here $\rho(q_i)$ is the unique extension
of the representation $\rho$ to affiliated elements.
\end{remark}

Another example of unbounded representations of $\Q_{2,\perp}(*)$ can be 
derived from 
Proposition~\ref{wild1}. Namely, let $\alpha_n\in{\mathbb R}$ be unbounded
sequence of numbers and let 
$$Q_1=\oplus_{n=1}^{\infty}\left(
\begin{array}{ccc}
1&1&\alpha_n\\
0&0&0\\
0&0&0
\end{array}\right),\  Q_2=\oplus_{n=1}^{\infty}\left(
\begin{array}{ccc}
0&-1&-1\\
0&1&1\\
0&0&0
\end{array}\right)$$ be operators on the Hilbert space
$H=\oplus_{n=1}^{\infty}H_n$, $H_n={\mathbb C}^3$.
$Q_1$, $Q_2$ determine a $*$-representation of $Q_{2,\perp}(*)$ defined on 
the set of, for example, finite vectors, i.e., vectors which are finite
linear combinations of vectors from $H_n$.

The similar example of representation can be constructed for the algebra
$Q_{3,3/2}(*)$.

We are not going to describe unbounded representations of
$Q_{2,\perp}(*)$, $\Q_{3,3/2}(*)$,
because  already the problem of
describing all their bounded $*$-representations is very complicated.

\section{Representations of algeb\-ras ${\mathcal Q}_{4,\lambda}$ and 
$*$-algebras
${\mathcal Q}_{4,\lambda}(*)$}

\subsection{Algebras ${\mathcal Q}_{4,\lambda}$ and their representations}
For each $\lambda\in{\mathbb C}$ the algebra ${\mathcal Q}_{4,\lambda}$ is 
non-zero. To see this we give a concrete example 
of idempotents $q_1,q_2,q_3,q_4$ as  operators defined on a linear space $X$.
The construction is a generalization of an example  given in \cite{bes}. 

Let $X$ be a linear space of complex-valued functions defined on ${\mathbb C}$.
Consider the operators $q_1$, $q_2$, $q_3$, $q_4\in L(X)$ defined by
\begin{eqnarray*}
&&(q_1f)(z)=z(f(z)+f(1-z)),\\
&&(q_2f)(z)=z(f(z)-f(1-z)),\\
&&(q_3f)(z)=(\lambda/2-z)f(z)+(1-\lambda/2+z)f(\lambda-1-z),\\
&&(q_4f)(z)=(\lambda/2-z)f(z)-(1-\lambda/2+z)f(\lambda-1-z).
\end{eqnarray*}
Simple computation shows that $q_1+q_2+q_3+q_4=\lambda I$ and $q_i^2=q_i$, 
$i=1,2,3,4$. This representations is  infinite dimensional.

Each algebra ${\mathcal Q}_{4,\lambda}$ is infinite dimensional. A linear 
basis for ${\mathcal Q}_{4,\lambda}$ is constructed in \cite{rss}.
It was also proved that for $\lambda\ne 2$ the algebras are not algebras
with the standard polynomial identities (PI-algebras), but for any
$x_1,\ldots,x_4\in\Q_{4,2}$ the following equality holds
$$\sum_{\sigma\in S_4}(-1)^{p(\sigma)}x_{\sigma(1)}x_{\sigma(2)}x_{\sigma(3)}
x_{\sigma(4)}=0$$
where $p(\sigma)$ is the parity of permutation $\sigma\in S_4$.
 
Let $\Lambda_{n,fd}$ be the set of all $\lambda\in{\mathbb C}$ for which 
there exists a finite-dimensional representation of $\Q_{n,\lambda}$.
The set $\Lambda_{4,fd}$ does not cover the whole complex plane, because
$Tr Q\in{\mathbb N}$ for any idempotent $Q$ in finite-dimensional space and therefore
$\Lambda_{n,fd}\subset{\mathbb Q}$ for any $n\in{\mathbb N}$. 
We have the following 
\begin{proposition}\label{q-4,alg}
$\Lambda_{4,fd}=\Lambda_{4,bd}=
\{2\pm\frac{2}{k},2\mid k\in{\mathbb N}\}$.
\end{proposition}
 \begin{proof}
Direct computation shows that
\begin{eqnarray*}
&\Q_{4,\lambda}={\mathbb C}\langle r_1,r_2,r_3,r_4,e\mid
r_k^2=e; \sum_{k=1}^4r_k=(2-\lambda)e\rangle=\\
&={\mathbb C}\langle x_1,x_2,x_3\mid\{x_1,x_2\}=x_3,\{x_1,x_3\}=x_2,
\{x_2,x_3\}=x_1,\\
& (\lambda-2)^2(x_1^2+x_2^2+x_3^2+1/4 e)=1\rangle,
\end{eqnarray*} 
where
\begin{eqnarray}\label{pro_1}
\begin{array}{l}
r_1=2q_1-e=(2-\lambda)(-x_1+x_2+x_3+1/2e),\\
r_2=2q_2-e=(2-\lambda)(x_1-x_2+x_3+1/2e),\\
r_3=2q_3-e=(2-\lambda)( x_1+x_2-x_3+1/2e),\\
r_4=2q_4-e=(2-\lambda)(-x_1-x_2-x_3+1/2e).
\end{array}
\end{eqnarray} 
If there exists a representation, $\pi$, of ${\mathbb C}\langle x_1,x_2,x_3,e
\mid\{x_1,x_2\}=x_3,\{x_1,x_3\}=x_2,
\{x_2,x_3\}=x_1,
 (\lambda-2)^2(x_1^2+x_2^2+x_3^2+1/4 e)=1\rangle$ by bounded operators in a 
Hilbert space $H$, then the operators $\pi(x_i)\otimes\sigma_i$, $i=1,2,3$,
 with
$\sigma_1= 
\left(\begin{array}{cc}
0&i\\
i&0\end{array}\right)$, $\sigma_2=
\left(\begin{array}{cc}
0&1\\
-1&0\end{array}\right)$ and  $\sigma_3=
\left(\begin{array}{cc}
-i&0\\
0&i\end{array}\right)$ define a representation of the universal enveloping
algebra $U(sl(2,{\mathbb C}))$ of the Lie algebra $sl(2,{\mathbb C})$ with 
an extra condition on the Casimir operator. Namely
if we take a basis $X_1$, $X_2$, $X_3$ in $sl(2,{\mathbb C})$ with the 
Lie bracket defined as 
\begin{equation}\label{li}
[X_1,X_2]=X_3,\  [X_2,X_3]=X_1,\  [X_3,X_1]=X_2.
\end{equation}
and denote by $\Delta$ the Casimir operator $X_1^2+X_2^2+X_3^2$ in 
$U(sl(2,{\mathbb C}))$ then $\rho(X_i)=\pi(x_i)\otimes\sigma_i$ is a 
representation of $U(sl(2,{\mathbb C}))$ so that
$(\lambda-2)^2\rho(\Delta)=(\lambda-2)^2/4-I$.
In the Hilbert space $H\otimes{\mathbb C}^2$ one
can choose an equivalent
scalar product  so that the operators $\rho(X_i)$
are skew-selfadjoint and $\rho$ is a $*$-representation
of the corresponding
$*$-algebra defined by the condition $X_i^*=-X_i$, $i=1,2,3$ 
(a $*$-representation of the Lie algebra $su(2)$)
It is known that $*$-representations of the $*$-algebra
such that the image of the Casimir operator equals 
$(1/4-1/(\lambda-2)^2)$  exist if and only if
$\lambda\in\{2\pm\frac{2}{k},2\mid k\in{\mathbb N}\}$. This implies that
$\Lambda_{n,fd}\subset\Lambda_{n,bd}\subset \{2\pm\frac{2}{k},2\mid k\in{\mathbb N}\}$.

To see the other inclusion, note that for any $\lambda\in
\{2\pm\frac{2}{k},2\mid k\in{\mathbb N}\}$ there exists a finite-dimensional
representation, $\pi$, of $U(sl(2,{\mathbb C}))$. Then
$\rho(x_i)=\pi(X_i)\otimes\sigma_i$, $i=1,2,3$, define a finite-dimensional 
representation of
${\mathbb C}\langle x_1,x_2,x_3,e
\mid\{x_1,x_2\}=x_3,\{x_1,x_3\}=x_2,
\{x_2,x_3\}=x_1,
 (\lambda-2)^2(x_1^2+x_2^2+x_3^2+1/4 e)=1\rangle$ and therefore
$\Q_{4,\lambda}$. The proof is complete.
\end{proof}

\begin{proposition}
Any finite-dimensional representation of $\Q_{4,\lambda}$, 
$\lambda\in\Lambda_{4,bd}\setminus\{2\}$ is equivalent
to a $*$-representation of ${\mathcal P}_{4,\lambda}$.
\end{proposition}
\begin{proof}
Let $\pi$ be a finite-dimensional representation
of $\Q_{4,\lambda}$, $\lambda\in\Lambda_{4,bd}\setminus\{2\}$, on a vector 
space $V_{\pi}$. The procedure is to find a scalar product, 
$(\cdot,\cdot)_{\pi}$
(a non-degenerate positive-definite hermitian form) such that
$$(\pi(q_i)\varphi,\psi)_{\pi}=(\varphi,\pi(q_i)\psi)_{\pi}, 
\quad\forall\ \varphi,\psi\in V_{\pi}, \ i=1,2,3,4,$$
or, equivalently,
$$(\pi(x_i)\varphi,\psi)_{\pi}=(\varphi,\pi(x_i)\psi)_{\pi}, 
\quad\forall\ \varphi,\psi\in V_{\pi}, \ i=1,2,3,$$
for the generators $x_i$ defined by (\ref{pro_1}).
 
Having the representation $\pi$ we can define a representation
$\rho$ of $U(sl(2,{\mathbb C})$ on $H(\rho)=V_{\pi}\otimes{\mathbb C}^2$ by 
setting $\rho(X_i)=\pi(x_i)\otimes\sigma_i$, $i=1,2,3$ for the generators $X_i$
of $U(sl(2,{\mathbb C})$ satisfying relations (\ref{li}).
It is known that $\rho$ is 
unitarizable, i.e., there exists a scalar product $(\cdot,\cdot)_{\rho}$ on
$H(\rho)$ such that $(\rho(X_i)\varphi,\psi)_{\rho}= 
-(\varphi,\rho(X_i)\psi)_{\rho}$.

Define a new representation $\tilde\pi$ of $\Q_{4,\lambda}$ on the Hilbert
space $H(\rho)\otimes{\mathbb C}^2$ by 
$\tilde\pi(x_i)=\rho(X_i)\otimes\sigma_i$, $i=1,2,3$. The scalar product, 
$(\cdot,\cdot)_{\tilde\pi}$ on
$H(\rho)\otimes{\mathbb C^2}$ is defined on elementary tensors  as
$(\varphi_1\otimes\psi_1,\varphi_2\otimes\psi_2)_{\tilde\pi}=
(\varphi_1,\varphi_2)_{\rho}(\psi_1,\psi_2)_{{\mathbb C}^2}$, for any
$\varphi_i\in 
H(\rho)$, $\psi_i\in{\mathbb C}^2$, where
 $(\cdot,\cdot)_{{\mathbb C}^2}$ is the standard scalar product on 
${\mathbb C}^2$.
Let $e_1$, $e_2$ be the standard basis vectors $(1,0)$ and $(0,1)$ respectively
in ${\mathbb C}^2$. For $\varphi,\psi\in V_{\pi}$ define
$$(\varphi,\psi)_{\pi}=((\varphi\otimes e_1)\otimes e_2+
(\varphi\otimes e_2)\otimes e_1,(\psi\otimes e_1)\otimes e_2+
(\psi\otimes e_2)\otimes e_1)_{\tilde\pi}.$$
It is easy to see that $(\cdot,\cdot)_{\pi}$ is a scalar product on $V_{\pi}$.
Moreover, $(\pi(x_i)\varphi,\psi)=(\varphi,\pi(x_i)\psi)$ for any
$\varphi$, $\psi\in V_{\pi}$, $i=1,2,3$.
We restrict ourselves by showing the last formula for the generator $x_1$.
\begin{eqnarray*}
&(\pi(x_1)\varphi,\psi)_{\pi}=\\&=
((\pi(x_1)\varphi\otimes e_1)\otimes e_2+(\pi(x_1)\varphi\otimes e_2)
\otimes e_1,(\psi\otimes e_1)\otimes e_2+
(\psi\otimes e_2)\otimes e_1)_{\tilde\pi}=\\
&=-((\pi(x_1)\otimes \sigma_1)(\varphi\otimes e_2)\otimes\sigma_1 e_1+
(\pi(x_1)\otimes \sigma_1)(\varphi\otimes e_1)\otimes\sigma_1 e_2,\\
&(\psi\otimes e_1)\otimes e_2+(\psi\otimes e_2)\otimes e_1)_{\tilde\pi}=\\
&=-((\rho(X_1)\otimes\sigma_1)((\varphi\otimes e_1)\otimes e_2+
(\varphi\otimes e_2)\otimes e_1),(\psi\otimes e_1)\otimes e_2+
(\psi\otimes e_2)\otimes e_1)_{\tilde\pi}=\\
&=-((\varphi\otimes e_1)\otimes e_2+
(\varphi\otimes e_2)\otimes e_1,(\rho(X_1)\otimes\sigma_1)
((\psi\otimes e_1)\otimes e_2+
(\psi\otimes e_2)\otimes e_1))_{\tilde\pi}=\\
&=-((\varphi\otimes e_1)\otimes e_2+
(\varphi\otimes e_2)\otimes e_1,\\
&(\pi(x_1)\otimes \sigma_1)(\psi\otimes e_2)
\otimes\sigma_1 e_1+
(\pi(x_1)\otimes \sigma_1)(\psi\otimes e_1)\otimes\sigma_1 e_2)_{\tilde\pi}=\\
&=((\varphi\otimes e_1)\otimes e_2+
(\varphi\otimes e_2)\otimes e_1,(\pi(x_1)\psi\otimes e_1)\otimes e_2+(\pi(x_1)\psi\otimes e_2)
\otimes e_1)_{\tilde\pi}=\\
&=(\varphi,\pi(x_1)\psi)_{\pi}.
\end{eqnarray*}
The proof is complete.
\end{proof}
As it follows from the previous proposition
finite-dimensional
indecomposable representations of $\Q_{4,\lambda}$ coincide with irreducible
$*$-representations of the $*$-algebra ${\mathcal P}_{4,\lambda}$, 
$\lambda\ne 2$.
For $\lambda=1+2/(2k+1)$ indecomposable representation
of ${\mathcal Q}_{4,\lambda}$ is unique, up to equivalence, and acts in a 
$(2k+1)$-dimensional vector
space. If $\lambda=1+2/(2k+2)$, there are four non-equivalent
representations of  ${\mathcal Q}_{4,\lambda}$ acting on $(k+1)$-dimensional
space (see \cite{os}[Section~2.2.1]). The algebra $\Q_{4,2}$ is wild and the 
problem of describing of its
 indecomposable representations is very complicated, \cite{bondar}.

\subsection{$*$-Algebras ${\mathcal Q}_{4,\lambda}(*)$ and their 
representations by bounded operators}

\begin{proposition}\label{wild2}
For $\lambda\in\Lambda_{4,bd}$ the $*$-algebra
$Q_{4,\lambda}(*)$ is $*$-wild.
\end{proposition}
\begin{proof}
That $Q_{4,\lambda}(*)$ is $*$-wild for $\lambda=1,2$ follows from 
Proposition~\ref{wild1}. Assume now that $\lambda\in\Lambda_{4,bd}
\setminus\{1,2\}$. Define
$p=(q_1+q_2)/2$, $q=(q_3+q_4)/2$, $r=(q_1-q_2)/2$, 
$s=(q_3-q_4)/2$. Direct computation shows that 
$Q_{4,\lambda}(*)$ is generated by $p$, $q$, $r$, $s$ and their adjoint and
the relations 
\begin{eqnarray}\label{rel_1}
&\begin{array}{cc}
pr=r(1-p),& ps=s(1-q),\\
r^2=p(1-p), & s^2=q(1-q),
\end{array}\\
&p+q=\lambda/2e\nonumber
\end{eqnarray}
From the relation it follows that $ps=s(-1+\lambda-p)$ and
$s^2=(\lambda/2-p)(1-\lambda/2+p)$.
To prove the statement we will construct a $*$-homomorphism
$\psi:\Q_{4,\lambda}(*)\to M_n({\mathbb C})\otimes C^*({\mathcal F}_2)$ 
($\simeq M_n(C^*({\mathcal F}_2))$ for some
$n\in{\mathbb N}$ depending on $\lambda$.

1. Let $\lambda=2+1/(2l)$, $l>0$. 
Let $E_n=\underbrace{\left(\begin{array}{ccc} e&&0\\&\ddots&\\
0&&e
\end{array}\right)}_{n\text{ times}}$, where $e$ is the identity in 
$C^*({\mathcal F}_2)$.

 Let
$J_3=\left(\begin{array}{c}A_1\\A_2\\A_3
\end{array}\right)$,  
$A_1=\frac{1}{N}\left(\begin{array}{ccccc}
e&0&0&0&0\\
0&e&0&0&0\\
0&0&2e&0&0\\
0&0&0&3e&0
\end{array}\right)$, 
$A_2=\frac{1}{N}\left(\begin{array}{ccccc}
e&0&e&e&e\\
0&2e&e&u_1&0\\
0&0&e&0&u_2\end{array}\right)$,
$A_3=\sqrt{E_5-A_1^*A_1-A_2^*A_2}$, where $u_1$, $u_2$ are the free unitary
generators of $C^*({\mathcal F}_2)$,
$N$ is chosen so that $||A_1^*A_1+A_2^*A_2||_{M_5(C^*({\mathcal F}_2))}<1$
($||\cdot||_{M_5(C^*({\mathcal F}_2))}$ is the $C^*$-norm on 
$M_5(C^*({\mathcal F}_2))$.

\vspace{0.1cm}

Define a $*$-homomorphism 
$\psi: \Q_{4,\lambda}(*)\to M_{12\cdot 2l}(C^*({\mathcal F}_2))$ as 
follows:
\begin{eqnarray*}
&\psi(p)=diag(\lambda_0E_{12},\lambda_1E_{12},\ldots,\lambda_{2l-1}E_{12}),\\
&\psi(r)=diag(0\cdot E_{12}, R_1,R_2,\ldots,R_{l-1},J_3J_3^*-1/2),\\
&\psi(s)=diag (S_0,S_1,\ldots,S_{l-1}),
\end{eqnarray*}
where
$\displaystyle\lambda_{2k}=1-\frac{k}{2l}$, 
$\lambda_{2k-1}=\displaystyle\frac{k}{2l}$, 
$R_k=\left(\begin{array}{cc}
0&\displaystyle \lambda_{2k-1}\lambda_{2k}E_{12}\\
E_{12}&0
\end{array}\right)$,
$$S_k=\left(\begin{array}{cc}
0&\displaystyle(\frac{\lambda}{2}-\lambda_{2k})(1-\frac{\lambda}{2}
+\lambda_{2k})E_{12}\\
E_{12}&0
\end{array}\right)
\
(k\ne l-1), \ S_{l-1}=\left(\begin{array}{cc}
0&s^{l-1}\\
s_{l-1}&0
\end{array}\right),$$
$$s_{l-1}=diag(x_1E_4,x_2E_3,x_3E_5),\quad
s^{l-1}=diag(y_1E_4,y_2E_3,y_3E_5), 
$$
$x_i$, $y_i$ are real numbers such that $x_i\ne x_j$, $y_i\ne y_j$ for
 $i\ne j$ and
$x_iy_i=(4l^2-1)/16l^2$. 

It is a routine to check the functor $F_{\psi}$ is full and we leave it to 
the reader.
The construction is similar for $\lambda= 2-1/2l$, $l\in{\mathbb N}$:
there exists a $*$-homomorphism $\psi:\Q_{4,\lambda}(*)\to M_{12\cdot 2l}
(C^*({\mathcal F}_2))$ with 
$\psi(p)=diag(\lambda_0E_{12},\lambda_1E_{12},\ldots,\lambda_{2l-1}E_{12})$, 
where $\lambda_{2k}=\displaystyle\frac{k}{2l}$,
$\lambda_{2k-1}=1-\displaystyle\frac{k}{2l}$.

2. Let $\lambda=2+2/(2l+1)$, $l>0$. 
Define $\psi: \Q_{4,\lambda}(*)\to M_{12(2l+1)}(C^*({\mathcal F}_2))$ as 
follows:

\begin{eqnarray*}
&\psi(p)=diag(\lambda_1E_{12},\lambda_2E_{12},\ldots,\lambda_{2l+1}E_{12}),\\
&\psi(r)=diag(0\cdot E_{12}, R_1,R_2,\ldots,R_{l-1},R_l),\\
&\psi(s)=diag (S_1,S_2,\ldots,S_{l},0\cdot E_{12}),
\end{eqnarray*}
where
$\displaystyle\lambda_{2k}=\frac{2k}{2l+1}$, 
$\displaystyle\lambda_{2k+1}=1-\frac{2k}{2l+1}$, $R_k=\left(\begin{array}{cc}
0&\displaystyle\lambda_{2k}\lambda_{2k+1}E_{12}\\
E_{12}&0
\end{array}\right)$,
 $(k\ne 1)$,
$$R_1=\left(\begin{array}{cc}
0&(2J_3J_3^*-1)\displaystyle\sqrt{\lambda_2\lambda_3}\\
(2J_3J_3^*-1)\displaystyle\sqrt{\lambda_2\lambda_3}&0
\end{array}\right),$$
$$S_k=\left(\begin{array}{cc}
0&\displaystyle(\frac{\lambda}{2}-\lambda_{2k-1})(1-\frac{\lambda}{2}+
\lambda_{2k-1})E_{12}\\
E_{12}&0
\end{array}\right),\ 
(k\ne 1,2), \ S_{k}=\left(\begin{array}{cc}
0&s^{k}\\
s_{k}&0
\end{array}\right), \ (k=1,2),$$ with
$$s_{k}=diag(x_1^kE_4,x_2^kE_3,x_3^kE_5),\quad
s^{k}=diag(y_1^kE_4,y_2^kE_3,y_3^kE_5)$$
$x_i^k$, $y_i^k$ are real numbers such that $x_i^k\ne x_j^k$, $y_i^k\ne y_j^k$
 for
 $i\ne j$ and
$x_i^ky_i^k=(\frac{\lambda}{2}-\lambda_{2k-1})(1-\frac{\lambda}{2}+
\lambda_{2k-1})$. 

One can check that the functor $F_{\psi}$ is full.
The construction is similar for $\lambda=2-2/(2l+1)$,  
$l\in{\mathbb N}$. In this case there exists a $*$-homomorphism
$\psi:\Q_{4,\lambda}\to M_{12(2l+1)}(C^*({\mathcal F}_2))$ with
$\psi(p)=diag(\lambda_1E_{12},\lambda_2E_{12},\ldots,\lambda_{2l+1}E_{12})$,
where $\lambda_{2k}=\displaystyle 1-\frac{2k}{2l+1}$, 
$\lambda_{2k+1}=\displaystyle
\frac{2k}{2l+1}$.

3. Let $\lambda=2+1/(2l+1)$, $l>0$. 
Define $\psi: \Q_{4,\lambda}(*)\to M_{12(2l+1)}(C^*({\mathcal F}_2))$ as 
follows:

\begin{eqnarray*}
&\psi(p)=diag(\lambda_1E_{12},\lambda_2E_{12},\ldots,\lambda_{2l+1}E_{12}),\\
&\psi(r)=diag(0\cdot E_{12}, R_1,R_2,\ldots,R_{l-1},R_l),\\
&\psi(s)=diag (S_1,S_2,\ldots,S_{l},J_3J_3^*-1/2),
\end{eqnarray*}
where
$\displaystyle\lambda_{2k}=\frac{k}{2l+1}$, 
$\displaystyle\lambda_{2k+1}=1-\frac{k}{2l+1}$, 
$R_k=\left(\begin{array}{cc}
0&\displaystyle\lambda_{2k}\lambda_{2k+1}E_{12}\\
E_{12}&0
\end{array}\right)$, $(k\ne l)$, $R_l=\left(\begin{array}{cc}
0&r^l\\
r_l&0
\end{array}\right)$ with 
$r_l=diag(x_1E_4,x_2E_3,x_3E_5),\quad
r^{l}=diag(y_1E_4,y_2E_3,y_3E_5)$, where $x_i$, $y_i$ are real numbers such that 
$x_i\ne x_j$, $y_i\ne y_j$
 for
 $i\ne j$ and
$x_iy_i=\frac{l}{2l+1}(1-\frac{l}{2l+1})$. 

$$S_k=\left(\begin{array}{cc}
0&\displaystyle(\frac{\lambda}{2}-\lambda_{2k-1})
(1-\frac{\lambda}{2} +\lambda_{2k-1})E_{12}\\
E_{12}&0
\end{array}\right)$$

The functor $F_{\psi}$ is full.
The construction is similar for $\lambda=2-2/(2l+1)$,  
$l\in{\mathbb N}$. In this case there exists a $*$-homomorphism
$\psi:\Q_{4,\lambda}\to M_{12(2l+1)}(C^*({\mathcal F}_2))$ with
$\psi(p)=diag(\lambda_1E_{12},\lambda_2E_{12},\ldots,\lambda_{2l+1}E_{12})$,
where
$\lambda_{2k}=\displaystyle 1-\frac{k}{2l+1}$, $\lambda_{2k+1}=\displaystyle
\frac{k}{2l+1}$.

\end{proof}

The problem of classification of all representations of $\Q_{4,\lambda}(*)$,
$\lambda\in\Lambda_{4,bd}$ is very difficult as it follows from 
Proposition~\ref{wild2}. However, if we restrict ourself to representations
such that the images of the generators $q_i$, $i=1, 2,3,4$ are selfadjoint, 
this problem reduces to the problem of describing
representations
of the $*$-algebra ${\mathbb C}\langle x_1,x_2,x_3,e
\mid\{x_1,x_2\}=x_3,\{x_1,x_3\}=x_2,
\{x_2,x_3\}=x_1,
 (\lambda-2)^2(x_1^2+x_2^2+x_3^2+1/4 e)=1, x_i^*=x_i\rangle$, $\lambda\ne 2$,
a factor $*$-algebra of the graded analogue of the Lie algebra
$so(3)$ (see section~4.1). Representations of the graded $so(3)$ are classified
in \cite{gorp}, see also \cite{os}. According to this result
there exists a unique, up to unitary
equivalence, irreducible representation of ${\mathcal P}_{4,\lambda}$ for
 $\lambda=1+2/(2k+1)$, acting  in a 
$(2k+1)$-dimensional vector
space and  there are four non-equivalent
representations of  ${\mathcal P}_{4,\lambda}$ acting on $(k+1)$-dimensional
space if $\lambda=1+2/(2k+2)$.
If $\lambda=2$, ${\mathcal P}_{4,2}$ has uncountable set of irreducible
unitarily non-equivalent representations which are one or two-dimensional
(see \cite{os}[Section~2.2.1]).

\subsection{Representations of ${\mathcal Q}_{4,\lambda}(*)$ by
unbounded operators}
As we already know, for each $\lambda\in{\mathbb C}$ the algebra
$\Q_{4,\lambda}$ and therefore $\Q_{4,\lambda}(*)$ is non-zero.
Representations of $\Q_{4,\lambda}$ or $\Q_{4,\lambda}(*)$ by bounded
operators on a Hilbert space  exist, however, not for all 
$\lambda\in{\mathbb C}$ 
(see Proposition~\ref{q-4,alg}).
\begin{proposition}\label{lambda4unbd}
$\Lambda_{4,unbd}={\mathbb C}$
\end{proposition}
\begin{proof}
In order to prove the statement it is enough for each $\lambda\in{\mathbb C}$
to give  a concrete construction
of unbounded representation of $\Q_{4,\lambda}(*)$. We follow \cite{rss}.

 Let $\varphi(\cdot)$, $\psi(\cdot)$ be the following idempotent 
matrix-functions
from  ${\mathbb C}$ to $M_2({\mathbb C})$:
$$\varphi(t)=\left(\begin{array}{cc}
t&t-t^2\\
1&1-t
\end{array}\right),
\psi(t)=\left(\begin{array}{cc}
t&-(t-t^2)\\
-1&1-t
\end{array}\right).$$
Consider a sequence of complex numbers $x_j=j(\lambda/2-1)$,
$j\in{\mathbb N}$.
Let $H=l_2$ and fix an orthonormal basis, $\{e_i,i\in{\mathbb N}\}$ in $H$.
Define operators $Q_i$, $Q_i^+$ on the set, $\Phi$, of finite vectors, 
i.e.\ finite
linear combinations of $e_i$, so that their matrix representations with 
respect to this fixed basis are given by
\begin{eqnarray*}
&Q_1=diag\{\varphi(x_1),\varphi(x_3),\varphi(x_5),\ldots\},\\
&Q_2=diag\{\psi(x_1),\psi(x_3),\psi(x_5),\ldots\},\\
&Q_3=diag\{1,\varphi(x_2),\varphi(x_4),\varphi(x_6),\ldots\},\\
&Q_4=diag\{1,\psi(x_2),\psi(x_4),\psi(x_6),\ldots\},\\
\end{eqnarray*}
 \begin{eqnarray*}
&Q_1^+=diag\{\varphi(x_1)^*,\varphi(x_3)^*,\varphi(x_5)^*,\ldots\},\\
&Q_2^+=diag\{\psi(x_1)^*,\psi(x_3)^*,\psi(x_5)^*,\ldots\},\\
&Q_3^+=diag\{1,\varphi(x_2)^*,\varphi(x_4)^*,\varphi(x_6)^*,\ldots\},\\
&Q_4^+=diag\{1,\psi(x_2)^*,\psi(x_4)^*,\psi(x_6)^*,\ldots\},\\
\end{eqnarray*}
where $A^*$ is the adjoint matrix to the matrix $A$.
Clearly, $\Phi$ is invariant with respect to $Q_i$, $Q_i^+$, $i=1,2,3,4$.
Moreover, direct calculation shows that $Q_i^2=Q_i$, $(Q_i^{+})^2=Q_i^+$,
$\sum_{i=1}^4 Q_i=\lambda I$, $\sum_{i=1}^4 Q_i^+=\overline{\lambda} I$.
Setting $\pi(q_i)=Q_i$, $\pi(q_i^*)=Q_i^+$, $D(\pi)=\Phi$, and then extending
$\pi$  to the whole algebra $\Q_{4,\lambda}(*)$ we obtain a $*$-representation
of $\Q_{4,\lambda}(*)$.
\end{proof}

We remark that the construction given in the proof can be  derived from 
one given in section~4.1.

In the same way as we obtained, in sections 4.1, 4.2, some
bounded representations
of $\Q_{4,\lambda}$ and $\Q_{4,\lambda}(*)$ from representations and,
respectively, unitary representations of the compact group $SU(2)$ (or
representations of the corresponding Lie algebra), unbounded representations of
$\Q_{4,\lambda}$ and $\Q_{4,\lambda}(*)$ can be obtained from representations
and, respectively, unitary representations of the Lie group 
$SL(2,{\mathbb R})$.

Let $\lambda\ne 2$ and let  $U(sl(2,{\mathbb C}))$ be the universal
enveloping  algebra of $sl(2,{\mathbb C})$ with the basis $X_1$, $X_2$, $X_3$
and the relations 
\begin{equation}\label{lie}
[X_1,X_2]=X_3,\  [X_2,X_3]=X_1,\  [X_3,X_1]=X_2.
\end{equation}
 Denote by
$\Delta$ the Casimir operator $X_1^2+X_2^2+X_3^2$. In the algebra
$U(sl(2,{\mathbb C}))\otimes
M_2({\mathbb C})$ consider the elements $x_1=
X_1\otimes\left(\begin{array}{cc}
0&i\\
i&0\end{array}\right)$, $x_2=
X_2\otimes\left(\begin{array}{cc}
0&1\\
-1&0\end{array}\right)$, $x_3=
X_3\otimes\left(\begin{array}{cc}
-i&0\\
0&i\end{array}\right)$. It is easy to check that they satisfy the relations
$\{x_1,x_2\}=x_3$,
$\{x_2,x_3\}=x_1$, $\{x_3,x_1\}=x_2$. 
We set
\begin{eqnarray}\label{pro}
\begin{array}{l}
\displaystyle Q_1=\frac{\lambda-2}{2}(-x_1+x_2+x_3)+\frac{\lambda}{4},\\
\displaystyle Q_2=\frac{\lambda-2}{2}(x_1-x_2+x_3)+\frac{\lambda}{4},\\
\displaystyle Q_3=\frac{\lambda-2}{2}(x_1+x_2-x_3)+\frac{\lambda}{4},\\
\displaystyle Q_4=\frac{\lambda-2}{2}(-x_1-x_2-x_3)+\frac{\lambda}{4}.
\end{array}
\end{eqnarray} 
Then, $Q_1+Q_2+Q_3+Q_4=\lambda I$. Moreover, $Q_i$, $i=1,2,3,4$, are 
idempotents iff
$$\hat\Delta=x_1^2+x_2^2+x_3^2=
\frac{1}{(\lambda-2)^2}-\frac{1}{4}.$$
It follows from the representation theory for $sl(2,{\mathbb C})$ that for any
 $\lambda\in{\mathbb C}$ there exists a representation 
of $U(sl(2,{\mathbb C}))$ 
such that the range of the Casimir operator is 
$(\frac{1}{4}-\frac{1}{(\lambda-2)^2})I$ (see \cite{V}). Namely, 
let $\chi=(l,\varepsilon)$, where $l$ is a complex number and $\varepsilon\in\{0,1/2\}$.
With each  such pair $\chi$ we associate  a space
$${\mathcal D}_{\chi}=\{f\in C^{\infty}({\mathbb R})\mid \hat f(x)=|x|^{2l}
(sgn x)^{2\varepsilon} f(\frac{1}{x})\in C^{\infty}({\mathbb R})\}.$$ 
Consider now representations $T_{\chi}$ of $SL(2,{\mathbb R})$ on ${\mathcal D}_{\chi}$ given by 
$$T_{\chi}(g)f(x)=|\beta x+\delta|^{2l}sgn^{2\varepsilon}(\beta x+\delta)
f(\frac{\lambda x+\gamma}{\beta x+\delta}),$$
where $g=\left(\begin{array}{cc}
\lambda&\beta\\
\gamma&\delta\end{array}\right)\in SL(2,{\mathbb R})$ (see \cite{ggv,V}). 
The infinitesimal
 operators of these representations are 
\begin{eqnarray*}
&A_1=2lx+(1-x^2)\frac{d}{dx},\\ 
&A_2=-2lx+(1+x^2)\frac{d}{dx},\\
& A_3=2l-2x\frac{d}{dx}.
\end{eqnarray*}
They satisfy the relations
$A_1A_2-A_2A_1=-2A_3$, $A_2A_3-A_3A_2=-2A_1$, $A_3A_1-A_1A_2=2A_2$. Let
$X_1=iA_1/2$, $X_2=A_2/2$, $X_3=-iA_3/2$. Then $X_1$, $X_2$, $X_3$ satisfy 
(\ref{lie})
 on ${\mathcal D}_{\chi}$ with $\Delta=X_1^2+X_2^2+X_3^2=-l(l+1)$. 
Obviously, a solution $l\in{\mathbb C}$ of $1/(\lambda-2)^2-1/4=l(l+1)$ exists
for any $\lambda\in{\mathbb C}$. 
Therefore there exist linear operators 
$Q_1$, $Q_2$, $Q_3$,
$Q_4$  on ${\mathcal D}_{\chi}$ such that $Q_1+Q_2+Q_3+Q_4=\lambda I$ and 
$Q_i^2=Q_i$.
We have the following expressions for $Q_i$, $i=1,2,3,4$:
\begin{eqnarray*}
Q_1=\frac{\lambda-2}{4}\left(\begin{array}{cc}
-A_3&A_1+A_2\\
A_1-A_2&A_3
\end{array}\right)+\frac{\lambda}{4},\ 
Q_2=\frac{\lambda-2}{4}\left(\begin{array}{cc}
-A_3&-A_1-A_2\\
A_2-A_1&A_3
\end{array}\right)+\frac{\lambda}{4},\\
Q_3=\frac{\lambda-2}{4}\left(\begin{array}{cc}
A_3&A_2-A_1\\
-A_1-A_2&-A_3
\end{array}\right)+\frac{\lambda}{4},\ 
Q_4=\frac{\lambda-2}{4}\left(\begin{array}{cc}
A_3&A_1-A_2\\
A_1+A_2&-A_3
\end{array}\right)+\frac{\lambda}{4}.
\end{eqnarray*}

The operators $Q_i$, $i=1,2,3,4$ define a representation of 
$\Q_{4,\lambda}$ on $D_{\chi}\otimes {\mathbb C}^2$.

\vspace{0.1cm}

It is known that for some values of $\chi$  one can introduce a scalar 
product $(\cdot,\cdot)$ 
in $D_{\chi}$ which is invariant with 
respect  to the representation $T_{\chi}$, i.e., $(\varphi,\psi)=
(T_{\chi}(g)\varphi, T_{\chi}(g)\psi)$, $\varphi,\psi\in D_{\chi}$.
The completion of $D_{\chi}$ with respect to the norm 
$||\varphi||^2=(\varphi,\varphi)$ gives us a Hilbert space $H_{\chi}$ 
and the continuous
extension of $T_{\chi}$ to $H_{\chi}$ gives a unitary representation of 
$SL(2,{\mathbb R})$.
In this case the infinitesimal operators  $A^{\chi}$ of $T_{\chi}$ will be
skew-selfadjoint, i.e. $(A^{\chi})^*=-A^{\chi}$. 

Recall that 
$D_{\chi}$ possesses an 
invariant
scalar product in the following cases: 

a) $l=-1/2+i\rho$, $\rho\in{\mathbb R}$,
$\varepsilon\in\{0,1/2\}$ and $\lambda=2+is$, $s\in{\mathbb R}$, $s\ne 0$, 
the corresponding 
representation $T_{\chi}$ is called a representation of the principal series;

b) $-1<l<0$, $l\ne -1/2$, $\varepsilon=0$ and $\lambda\in(-\infty, 0)\cup 
(4,+\infty)$, the corresponding representation $T_{\chi}$ is called a 
representation of a supplementary series. 

For $l\in\frac{1}{2}{\mathbb Z}$, 
$l\leq -1$ and $\varepsilon$ satisfying the condition
 $l+\varepsilon\in{\mathbb Z}$ the corresponding space $D_{\chi}$ has  two
subspaces $F_l^+$, $F_l^-$ which are invariant with respect to the operators 
$T_{\chi}(g)$. It is not possible to introduce an  invariant scalar product
on $D_{\chi}$, but it is possible to do it on each of these subspaces.
The corresponding subrepresentations of $T_{\chi}$ are called representations
of the discrete series. In this case $\lambda$ takes values $2\pm \frac{2}{k}$,
$k\in{\mathbb N}$.

Taking the infinitesimal representation of a  unitary representation of 
$SL(2,{\mathbb R})$ in a Hilbert space $H$ we define the (unbounded) 
operators $Q_i$, $i=1,2,3,4$ on $H\oplus H$ as above. These operators are
densely defined and clearly, 
 there exist a dense invariant domain $D$ such that the 
 the equalities
$(Q_1+Q_2+Q_3+Q_4)\varphi=\lambda\varphi$, $Q_i^2\varphi=Q_i\varphi$ hold
for any $\varphi\in D$.
Moreover, if $\lambda$ is real then $Q_1^*\supseteq Q_3$, $Q_2^*\supseteq Q_4$.


\subsubsection{Unbounded idempotents the sum of which  is zero}

The rest of the section is devoted to a detailed discussion of unbounded
representations of $\Q_{4,0}(*)$.
Note that idempotents whose sum is zero were studied in \cite{bes} in 
connection with their investigation of the concept of logarithmic residues in 
Banach algebras.
We will study representations of $\Q_{4,0}(*)$ under some additional 
conditions which will allow us to classify them up to unitary equivalence. 


Having these unbounded representations we  construct a $C^*$-algebra ${\mathfrak A}$ and unbounded elements $q_1$,$q_2$, $q_3$, $q_4$ which are affiliated
with ${\mathfrak A}$ such that any non-degenerate representation of
${\mathfrak A}$ extended to its affiliated elements gives us an integrable
representation of $\Q_{4,0}(*)$ defined below. Moreover, any integrable
representation with some extra condition can be obtained this way. 


\vspace{0.1cm}

{\bf 1. Hilbert space level.}

\vspace{0.1cm}

Consider new generators in $\Q_{4,0}(*)$ given by
$p=(q_1+q_2)/2$, $q=(q_3+q_4)/2$, $r=(q_1-q_2)/2$, 
$s=(q_3-q_4)/2$ and their adjoint. Direct computation shows that 
relations $q_1+q_2+q_3+q_4=0$, $q_i^2=q_i$, $i=1,\ldots,4$, are 
equivalent
to the following ones: 
\begin{eqnarray}\label{rel1}
\begin{array}{cc}
pr=r(1-p),& ps=s(-1-p),\\
r^2=p(1-p), & s^2=-p(p+1).
\end{array}
\end{eqnarray}
Let ${\mathcal A}_{4,0}$ be the quotient of $\Q_{4,0}(*)$ by the two-sided
$*$-ideal generated by the elements $pr^*-rp$, $ps^*-sp$ and $p-p^*$.
So we have additional relations in ${\mathcal A}_{4,0}$, namely,
\begin{equation}\label{rel2}
pr^*=rp,\ ps^*=sp,\ p=p^*.
\end{equation}
 In what follows we will study  representations of ${\mathcal A}_{4,0}$. 
Obviously, any 
$*$-representation of ${\mathcal A}_{4,0}$ is a representation of 
$\Q_{4,0}(*)$.

Henceforth, $H_a(A_1,\ldots,A_n)$ will denote the set of
joint analytic vectors for selfadjoint operators $A_1,\ldots,A_n$ (see 
\cite{Shbook}). Also, a 
linear set $\Phi$ will
be called  a core for a closed operator $A$ if $\Phi\subseteq D(A)$ and
the closure of the operator $A$ restricted to the domain $\Phi$ is equal 
$A$. 
\begin{definition}\label{def}
We say that closed operators ($p=p^*$, $q=q^*$, $r$, $s$, $r^*$, $s^*$) is a 
representation of commutation relations $(\ref{rel1})-(\ref{rel2})$
on a Hilbert space $H$ if there exists a linear dense subset $\Phi\subset
H$ such that

$1$. $\Phi\subseteq H_a(p, q, r^*r, s^*s)$;

$2$. $\Phi$ is a core for the operators $r$, $r^*$, $s$ and $s^*$;

$3$. relations $(\ref{rel1})-(\ref{rel2})$ hold on $\Phi$.
\end{definition}
A family $\{A_j \mid j \in J\}$ of closed unbounded operators on a Hilbert space
$H$ is called {\it irreducible} if decomposition $A_j=B_j \oplus C_j$ for all
$j\in J$ with respect to an orthogonal direct sum $H=H_1 \oplus H_2$ is only possible
when either $H_1 = \{0\}$ or $H_2 = \{0\}$, or equivalently if
$$\{C\in B(H)\mid C A_j
\subseteq A_j C \ 
\mbox{and}\ C^* A_j \subseteq A_j C^*, j \in J\}=
{\mathbb C}I.$$
These and other definitions and facts from the general theory
of unbounded representations of algebras and relations can be found, 
for example,
in \cite{Shbook}.
\begin{remark}\rm
The operators $p=p^*$, $q=q^*$, $r$, $s$, $r^*$, $s^*$  satisfying
the conditions of Definition~\ref{def} define a 
representation,
$\pi$, of ${\mathcal A}_{4,0}$ and $\Q_{4,0}(*)$ on the domain $\Phi$ in the 
sense of definition given section~2.3.
This domain is not unique and therefore there are many $*$-representations
of ${\mathcal A}_{4,0}$ corresponding to the closed operators 
$p=p^*$, $q=q^*$, $r$, $s$, $r^*$, $s^*$. Among them there is a unique 
selfadjoint representation $\pi^*$.
This representation is irreducible iff the family of the closed operators
is irreducible, two such representations are unitarily equivalent iff
the corresponding families of closed operators are unitarily equivalent
(see Remark~\ref{rep_1}).

In what follows we mean these selfadjoint representations when we talk about 
integrable representation of the $*$-algebra ${\mathcal A}_{4,0}$.
\end{remark}

Let $O_x$ be the trajectory of the point $x$ with respect 
to the mappings $F_1(x)=1-x$, $F_2(x)=-1-x$, i.e.,
$O_{x}=\{F_{i_1}\ldots F_{i_n}(x)\mid i_k\in\{1,2\}, 
n\in{\mathbb N}\}=\{(-1)^n(x-n), (-1)^n(-x-n)\mid 
n\in{\mathbb N}\cup\{0\}\}$.

Let $O_{0}^+=\{F_{i_1}\ldots F_{1}(0)\mid i_k\in\{1,2\}, 
n\in{\mathbb N}\}=\{(-1)^{k+1}k\mid k\in{\mathbb N}\}$ and
$O_{0}^-=\{F_{i_1}\ldots F_{2}(0)\mid i_k\in\{1,2\}, 
n\in{\mathbb N}\}=\{(-1)^{k}k\mid k\in{\mathbb N}\}$. 

Denote by $l_2(K)$ the separable Hilbert space 
with the orthonormal basis $\{e_{\mu}\}_{\mu\in K}$

\begin{theorem}\label{rep}
Any irreducible integrable representation $\pi$ of the $*$-algebra 
${\mathcal A}_{4,0}$ in a 
Hilbert 
space $H$ such that $\ker p\ne\{0\}$ is unitarily equivalent to one of the 
following:

I. $H=l_2(O_{\lambda})$ 
\begin{equation}\label{I}
\begin{array}{lll}
pe_{\mu}&=&\mu e_{\mu}\\
qe_{\mu}&=&-\mu e_{\mu}\\
re_{\mu}&=&(1-\mu)e_{1-\mu}\\
se_{\mu}&=&-(1+\mu)e_{-1-\mu}\end{array}\quad
\end{equation}
where $\lambda\in(-1/2,1/2)\setminus\{0\}$.

II.
$H=l_2(O_{1/2})$
\begin{equation}\label{II}
\begin{array}{lll}
pe_{\mu}&=&\mu e_{\mu}\\
qe_{\mu}&=&-\mu e_{\mu}\\
re_{\mu}&=&\left\{
\begin{array}{ll}
ae_{1/2}&\mu=1/2\\
(1-\mu)e_{1-\mu}& \mu\ne 1/2
\end{array}\right.\\
se_{\mu}&=&-(1+\mu)e_{-1-\mu}\end{array}\quad
\end{equation}
where $a=\pm 1/2$, $s=\pm 1$.

III.
$H=l_2(O_{-1/2})$
\begin{equation}\label{III}
\begin{array}{lll}
pe_{\mu}&=&\mu e_{\mu}\\
qe_{\mu}&=&-\mu e_{\mu}\\
re_{\mu}&=&(1-\mu)e_{1-\mu}\\
se_{\mu}&=&\left\{
\begin{array}{ll}
ae_{-1/2}&\mu=-1/2\\
-(1+\mu)e_{-1-\mu}& \mu\ne -1/2
\end{array}\right.\end{array}\quad
\end{equation}
where $a=\pm 1/2$, $s=\pm 1$.

IV. $H=l_2(O_{0}^-)$
\begin{equation}\label{IV}
\begin{array}{lll}
pe_{\mu}&=&\mu e_{\mu}\\
qe_{\mu}&=&-\mu e_{\mu}\\
re_{\mu}&=&(1-\mu)e_{1-\mu}\\
se_{\mu}&=&\left\{
\begin{array}{ll}
0&\mu=-1\\
-(1+\mu)e_{-1-\mu}& \mu\ne -1
\end{array}\right.\\
\end{array}\quad
\end{equation}

V. $H=l_2(O_{0}^+)$
\begin{equation}\label{V}
\begin{array}{lll}
pe_{\mu}&=&\mu e_{\mu}\\
qe_{\mu}&=&-\mu e_{\mu}\\
re_{\mu}&=&\left\{
\begin{array}{ll}
0&\mu=1\\
(1-\mu)e_{1-\mu}& \mu\ne 1
\end{array}\right.\\
se_{\mu}&=&-(1+\mu)e_{-1-\mu}
\end{array}\quad
\end{equation}

\end{theorem}
\begin{proof}
Let $p$, $q$, $r$, $s$ be closed operators satisfying the conditions 
$(1)-(3)$ of  Definition~\ref{def}, $E_p(\cdot)$  the resolution
of the identity for the selfadjoint operator $p$ and $r=u_r|r|$, $s=u_s|s|$ 
the polar decompositions of the closed operators $r$, $s$. Here 
$|r|=(r^*r)^{1/2}$, $|s|=(s^*s)^{1/2}$, $\ker u_r=\ker r=\ker |r|$ and
$\ker u_s=\ker s=\ker |s|$. By \cite{os},  we conclude that
$$\mbox{ $p$, $|r|$ and $p$, $|s|$ commute strongly},$$
(i.e. in the sense of resolutions of the identities)
\begin{equation}{\label{res}}
E_p(\Delta)u_r=u_rE_p(1-\Delta),\ E_p(\Delta)u_s=u_sE_p(-1-\Delta),\ 
\forall\Delta\in{\mathfrak B}({\mathbb R}).
\end{equation}
Here ${\mathfrak B}({\mathbb R})$ is the Borel $\sigma$-algebra on 
${\mathbb R}$. Assume first that $\ker p(1-p)(1+p)=\{0\}$.
Since $r^2\varphi=p(1-p)\varphi$ and $s^2\varphi=-p(1+p)\varphi$, 
$\varphi\in\Phi$ we have that
$\ker r\subset\ker r^2=\{0\}$, $\ker s\subset\ker s^2=\{0\}$  and
 $u_r$, $u_s$ are unitary operators.
 The equality  $pr^*\varphi=rp\varphi$ gives  $pr^*r\varphi=r^2(1-p)\varphi$ 
 which implies $|r|=|1-p|$. From (\ref{res}) one can easily derive that
$u_r{\mathcal D}(|1-p|)\in {\mathcal D}(|p|)$ and 
$$|p|u_r\psi=u_r|1-p|\psi,\quad \psi\in {\mathcal D}(|1-p|)={\mathcal D}(p).$$
From $r^2\varphi=p(1-p)\varphi$ we have $u_r|r|r\varphi=p(1-p)\varphi$
for any $\varphi\in\Phi$ and, since $r\varphi\in {\mathcal D}(p)$, 
$|p|u_r^2|1-p|\varphi=p(1-p)\varphi$ and $u_r^2=\mbox{sgn}(p(1-p))$. 
Setting $u_1=\mbox{sgn}(p)u_r$, we get $r=u_1(1-p)$. Similarly,
$|s|=|1+p|$,  $u_s^2=\mbox{sgn}(-p(1+p))$ and $s=-u_2(1+p)$, where
$u_2=\mbox{sgn}(p)u_s$.

It follows from (\ref{res}) that if $\Delta\in {\mathfrak B}({\mathbb R})$ is 
invariant with respect to the mappings $F_1(\lambda)=1-\lambda$, 
$F_2(\lambda)=-1-\lambda$ then $E_p(\Delta)$ commutes with $p$, $r$, $s$,
$r^*$, $s^*$ in the sense $E_p(\Delta)T\subseteq TE_p(\Delta)$  for 
$T=p,q,r,s,r^*,s^*$. Therefore,  if ($p$, $q$, $r$, $s$, $r^*$, $s^*$) is 
irreducible then 
$E_p(\Delta)=cI$, where $c=0,1$. One can easily check that the set 
$\tau=[-1/2,1/2]$ 
intersects every trajectory $O_{\lambda}$ exactly in one points which implies
that the spectral measure $E_p(\cdot)$ is concentrated on an orbit 
$O_{\lambda}$ for some $\lambda\in[-1/2, 1/2]$ if the representation is 
irreducible.
 Searching irreducible representation we can assume now that the representation
space $H$ is
a direct sum $\oplus_{\mu\in O_{\lambda}}H_{\mu}$, where $H_{\mu}$ is 
an eigenspace of $p$ corresponding to the eigenvalue $\mu$.

If $\lambda\ne 0, \pm 1/2$  one can easily check
that the linear span of the vectors $\{u_{i_1}\ldots u_{i_k}e\mid i_l\in\{r,s\}, 
k\in {\mathbb N}\cup\{0\}\}$,
where $e\in H_{\lambda}$, 
is invariant with respect to the operators $p$, $q$, $r$, $s$, $r^*$, $s^*$ 
and moreover the operators restricted to the closure of this subspace define an
irreducible representation which is given by
formulae (\ref{I}). 

If $\lambda=1/2$ then $u_rH_{1/2}\subset H_{1/2}$. It is not difficult to see
that, given an irreducible representation, the operator $u_r$ has
an eigenvector $e \in H_{1/2}$ and the vectors  
$\{u_{i_1}\ldots u_{i_l}u_se\mid i_k\in{r,s},k\in{\mathbb N}\cup\{0\}\}$ 
build a basis of the representation space
The corresponding irreducible 
representation is given by (\ref{II}).

Representations related to the orbit $O_{-1/2}$ can be obtained in a similar 
way.
Note that there is no representation related with the trajectory $O_{0}$ 
such that $\ker p(1-p)(1+p)=\{0\}$. 

If $\ker(1-p)(1+p)\ne\{0\}$ and  $\ker p=\{0\}$ then using the same arguments 
one 
can show that $r^*r=|1-p|$ $s^*s=|1+p|$, $u_r^2=\mbox{sgn}(p(1-p))$, 
$u_s^2=\mbox{sgn}(-p(1+p))$ and
$u_r^*|_{ker(1-p)}=0$, $u_s^*|_{ker(1+p)}=0$. 
Moreover $W_1=\oplus_{k\in{\mathbb N}} u_s(u_ru_s)^k\ker(1-p)\oplus
\oplus_{k\in{\mathbb N}} (u_ru_s)^k\ker(1-p)$ and 
$W_2=\oplus_{k\in{\mathbb N}} u_r(u_su_r)^k\ker(1+p)\oplus
\oplus_{k\in{\mathbb N}} (u_su_r)^k\ker(1+p)$ are invariant with respect to
the operators $p$, $q$, $r$, $s$, $r^*$, $s^*$  and the corresponding 
irreducible representation are given by (\ref{V}) and (\ref{IV}) 
respectively.  

Clearly, $W_1\perp W_2$ and any representation space $H$ can be 
decomposed  into a direct sum of invariant with respect to the representation
subspaces, namely, 
 $H=W_1\oplus W_2\oplus W_3$, where $W_3=(W_1\oplus W_2)^{\perp}$.
Moreover, if $\ker p =\{0\}$, we obtain $\ker p(1-p)(1+p)|_{W_3}=\{0\}$.
Setting $u_1=\mbox{sgn}(p)u_r$ on $(\ker (1-p))^{\perp}$, 
$u_2=\mbox{sgn}(p)u_s$ on $(\ker (1+p))^{\perp}$ and extending them  to
$\ker (1-p)$ and $\ker (1+p)$ in a way that $u_1$, $u_2$ are unitary and 
satisfying (\ref{res}) we get that the operators $\tilde p=pp_1$, 
$\tilde q=qp_1$,
$\tilde r=u_1(1-p)p_1$, $\tilde s=-u_2(1+p)p_1$, where $p_1$ is the projection
onto $W_1$, define a  representation of ${\mathcal A}_{4,0}$ on $W_1$, and
$\hat p=pp_2$, $\hat q=qp_2$
$\hat r=u_1(1-p)p_2$, $\hat s=-u_2(1+p)p_2$, where $p_2$ is the projection 
onto $W_2$, define a representation of ${\mathcal A}_{4,0}$ on $W_2$.
Moreover, any representation on $W_1$ and $W_2$ can be obtained  this way.
The proof is finished.
\end{proof}

\vspace{0.1cm}

{\bf 2.  $C^*$-algebra level}

\vspace{0.1cm}

In the sequel, we use the following notation. The set of multiplier of a 
$C^*$-algebra $A$ is denoted by $M(A)$. The notation $T\eta A$ means $T$ is
affiliated with the algebra $A$ and $z_T$ denotes its $z$-transform. We write
$Mor(A,B)$ for the set of morphisms from $A$ to another $C^*$-algebra
$B$. For the definition  and facts related to these notions we 
refer the reader to \cite{wor1}.

It follows from the proof of Theorem~\ref{rep} that any representation 
($p$, $q$, $r$, $s$) in a Hilbert space $H$ provided $\ker p=\{0\}$ is of the 
form: 
$a=a^1\oplus a^2\oplus a^3$, where $a\in\{p,q,r,s, r^*, s^*\}$, $a^i$ are
operators on $W_i$ , $i=1,2,3$ described in the proof.
Moreover,
$p^i=(p^i)^*$, $q^i=-p^i$, $\ker (1-p^3)p^3(1+p^3)=\{0\}$, 
$Sp(p^1)\subset O_0^+$, $Sp(p^2)\subset O_0^-$, $r^i=u_1^i(1-p^i)$, 
$s^i=-u_2^i(1+p^i)$, where $u_1^i$, $u_2^i$ are unitary operators such that
\begin{equation}\label{cross}
(u_1^i)^2=1, \ (u_2^i)^2=1, \ (u_1^i)^*p^iu_1^i=1-p^i,\  (u_2^i)^*pu_2^i=-1-p^i,\ i=1,2,3.
\end{equation}
Our aim now is  to define a $C^*$-algebra ${\mathfrak A}$ generated by 
selfadjoint
element $p=p^*$ and unitary elements $u_1$, $u_2$ satisfying (\ref{cross}) 
and affiliated with ${\mathfrak A}$. Namely, we look for $C^*$-algebra with the 
following universal property: for any $C^*$-algebra ${\mathfrak A}'$ and any 
$U_1$, 
$U_2$, $P\ \eta\ {\mathfrak A}'$ such that $U_1$, $U_2$ 
are unitary, $P$ is selfadjoint  and $U_1^*PU_1=1-P$, $U_2^*PU_2=-1-P$,
$U_1^2=1$, $U_2^2=1$, there exists unique $\Phi\in Mor({\mathfrak A}, {\mathfrak A}')$ 
such that $\Phi(u_1)=U_1$, $\Phi(u_2)=U_2$, $\Phi(p)=P$.

Let ${\mathfrak A}=C_{\infty}({\mathbb R})\ltimes_{\alpha}
({\mathbb Z}_2\times{\mathbb Z}_2)$, where $C_{\infty}({\mathbb R})$ is the 
algebra of all continuous, vanishing at infinity functions on ${\mathbb R}$.
The action $\alpha$ of ${\mathbb Z}_2\times{\mathbb Z}_2$ on 
$C_{\infty}({\mathbb R})$ is defined by 
$$(gf)(\xi)=f(g\xi),$$
where the action on the real line ${\mathbb R}$ is given by
$$g_1\xi=1-\xi,\ g_2\xi=-1-\xi$$
for the generators $g_1\in{\mathbb Z}_2$, $g_2\in{\mathbb Z}_2$.

Then there exist unitary operators $u_1$, $u_2\in M({\mathfrak A})$ such that
$$u_1^*fu_1=g_1f, \ u_2^*fu_2=g_2f, u_1^2=1,\ u_2^2=1.$$
Let now $p$ be the function defined by $p(\xi)=\xi$ for all 
$\xi\in{\mathbb R}$. Clearly, $p\ \eta\  C_{\infty}({\mathbb R})$, 
$(1-p)\ \eta\  C_{\infty}({\mathbb R})$, $(-1-p)\ \eta\ C_{\infty}({\mathbb R})$ and since the 
embedding $C_{\infty}({\mathbb R})\hookrightarrow 
C_{\infty}({\mathbb R})\ltimes_{\alpha}
({\mathbb Z}_2\times{\mathbb Z}_2$) is in $Mor(C_{\infty}({\mathbb R}),{\mathfrak A})$ 
we have $p,(1-p), (-1-p)\ \eta\ {\mathfrak A}$ and $u_1^*pu_1=g_1p$, $u_2^*pu_2=g_2p$.

Clearly, ${\mathfrak A}$ possesses the universality property defined above.

\begin{proposition}
The elements $r=u_1(1-p)$, $s=-u_2(1+p)$ are affiliated with ${\mathfrak A}$. 
Moreover, there is a dense domain, $D$, of ${\mathfrak A}$ such that  
relations (\ref{rel1})  hold on $D$.
\end{proposition}
\begin{proof}The first statement follows from \cite[Example~2]{wor1} and 
the fact that
$u_1$, $u_2\in M({\mathfrak A})$ are invertible and $(1-p), (1+p)\ \eta\ 
{\mathfrak A}$. In this case $D(r)=D(s)=D(p)$. One can easily check also that
the relations hold on $D(p^2)$ which is dense in ${\mathfrak A}$.
\end{proof}

By \cite{wor2}[Theorem~3.3], the affiliated elements $p$, $q$, $r$, $s$
generate the $C^*$-algebra ${\mathfrak A}$: for any Hilbert space
$H$, any $C^*$-subalgebra, $B$, of $B(H)$ and any non-degenerate 
representation
$\pi$ of ${\mathfrak A}$ on $H$ we have $\pi(X_i)\ \eta\ B$ for
$X_i=p,q,r,s$ implies $\pi\in Mor({\mathfrak A},B)$. 

We see also that any representation of ${\mathfrak A}$ generates a 
representation
of ${\mathcal A}_{4,0}$ satisfying the conditions of Definition~\ref{def}. 
Moreover, any such irreducible representation is unitarily equivalent either to
one from Theorem~\ref{rep} or to one-dimensional zero representation 
$\pi(x)=0$, $x=p$, $q$, $r$, $s$. Conversely, for any representation
$P$, $Q$, $R$, $S$, $R^*$, $S^*$ of ${\mathcal A}_{4,0}$ defined in 
Definition~\ref{def} and such that $\ker P=\{0\}$ there exists a 
representation $\pi$ of ${\mathfrak A}$ having the property $X=\pi(x)$,
$(X,x)=(P,p),(Q,q),(R,r),(S,s) (R^*,r^*), (S^*,s^*)$, where $\pi(x)$ is the 
unique extension of $\pi$ to the affiliated elements.

One can also define idempotents $q_1$, $q_2$, $q_3$, $q_4$ with the  zero 
sum 
 in a way that all of them are affiliated with ${\mathfrak A}$.
 Further we will use the following statement from \cite{wor1}:

Let $A$ be a $C^*$-algebra; $a$, $b$, $c$, $d\in M(A)$ and 
$Q=\left(
\begin{array}{cc}
d,&-c^*\\
b&a^*
\end{array}\right)$. Assume that 
(1) $ab=cd$, (2)
$a^*A$ is dense in $A$, (3) $dA$ is dense in $A$, (4)
$Q(A\oplus A)$ is dense in $A\oplus A$.
Then there exists $T\ \eta\ A$ such that
1. $dA$ is a core for $T$ and 
$Tdx=bx$
for any $x\in A$.
2. For any $x,y\in A$
$$
\left(\begin{array}{c}
x\in D(T)\ \mbox{and}\\
y=Tx
\end{array}\right)\Leftrightarrow (ay=cx)
$$
If $Q$ is invertible then $D(T)=dA$.

Using this statement we prove  the following 
\begin{proposition}
There exists elements $q_1$, $q_2$, $q_3$, $q_4\ \eta\ {\mathfrak A}$ and a dense
domain $D$ such that $D$ is invariant with respect to $q_i$, $q_i^*$,
$D$ is a core for any $q_i$ and $q_i^*$ and  
the relations $\sum_{i=1}^4q_i=0$, $\sum_{i=1}^4q_i^*=0$ and
$q_i^2=q_i$, $(q_i^*)^2=q_i^*$ hold on $D$. Moreover, $q_1x=px+rx$, 
$q_2x=px-rx$, 
$q_3x=-px+sx$, $q_4x=-px-sx$
for any $x\in D$.
\end{proposition}
\begin{proof} First we will prove the existence of $q_1$. 
Let 
$$a=d=(1-z_p^2)^{1/2}(1-z^2_{1-p})^{1/2},\ 
b=c=z_p(1-z_{1-p}^2)^{1/2}+u_1z_{1-p}(1-z_p^2)^{1/2}$$
Then, clearly, $ab=cd$ and $a^*{\mathfrak A}=d{\mathfrak A}$ is dense in 
${\mathfrak A}$. We have also
$$Q^*Q=\left(
\begin{array}{cc}
d^*d+b^*b,&0\\
0,&cc^*+aa^*
\end{array}\right).$$ We state that $Q^*Q({\mathfrak A}\oplus{\mathfrak A})$ is
 dense in ${\mathfrak A}\oplus{\mathfrak A}$. It is 
enough to see that $(d^*d+b^*b){\mathfrak A}$ and $(cc^*+aa^*){\mathfrak A}$ 
are dense
in ${\mathfrak A}$. Assume that
$\overline{(d^*d+b^*b){\mathfrak A}}\ne{\mathfrak A}$. Then there exists a pure
state $w$ on ${\mathfrak A}$ such that $w((d^*d+b^*b)x)=0$ for any $x\in{\mathfrak A}$.
Let $\pi$ be the GNS representation of ${\mathfrak A}$ acting on a Hilbert space
$H_{\pi}$ and $\Omega\in H_{\pi}$ be the corresponding cyclic vector such that
$$w(x)=(\Omega,\pi(x)\Omega)$$
for any $x\in {\mathfrak A}$. This gives that the range 
$R(\pi(d)^*\pi(d)+\pi(b)\pi(b)^*)$ belongs to the set
$\{\varphi\in H_{\pi}\mid (\Omega,\varphi)=0\}$. Since the operators 
$\pi(d)^*\pi(d)$ and  $\pi(b)\pi(b)^*$ are positive and commute with each 
other, we have that $\Omega\in \ker \pi(b)\pi(b)^*\cap\ker \pi(d)^*\pi(d)$. 
This contradicts  the statement 
that
$d{\mathfrak A}$ is dense in ${\mathfrak A}$.
Using the same arguments one can show that 
$\overline{(cc^*+aa^*){\mathfrak A}}={\mathfrak A}$.
 The statement about density of $Q({\mathfrak A}\oplus {\mathfrak A})$ can be 
easily derived from the 
density of $Q^*Q(A\oplus A)$.

Let $q_1\ \eta\ {\mathfrak A}$ be the operator from the previous statement. 
Then 
$d{\mathfrak A}$ is a core for $q_1$ and 
$q_1x=px+rx$ for any $x\in d{\mathfrak A}$.   
Similarly, we can construct $q_i\ \eta\ {\mathfrak A}$, $i=2,3,4$ such that
$d{\mathfrak A}$ is a core for $q_2$, $q_2x=px-rx$ for any $x\in d{\mathfrak A}$ 
and $d'{\mathfrak A}:=(1-z_p^2)^{1/2}(1-z^2_{1+p})^{1/2}{\mathfrak A}$ is a core
 for $q_3$ and $q_4$, $q_3x=-px+sx$, $q_4x=-px-sx$ for any $x\in d'{\mathfrak A}$. 
Set $D=D(p^2)$. Then $D=d{\mathfrak A}=d'{\mathfrak A}$ and $D$ is a core for 
all idempotents $q_1$, $q_2$, $q_3$, $q_4$. Moreover, the relations 
$q_1+q_2+q_3+q_4=0$, $q_i^2=q_i$ hold on $D$.
\end{proof}

\begin{remark}\rm
It was proved in \cite{bes} that $\Q_{4,0}$ and therefore $\Q_{4,0}(*)$ is not 
trivially $B$-representable, i.e.,
there exists no non-trivial isomorphism of $\Q_{4,0}$ into a subalgebra
of a Banach algebra and respectively $*$-isomorphism of $\Q_{4,0}(*)$ into 
a $*$-subalgebra of an involutive  Banach algebra .
We have shown that there exist a $C^*$-algebra ${\mathfrak A}$ and unbounded
elements $q_1$, $q_2$, $q_3$, $q_4$ which are affiliated with ${\mathfrak A}$
and such that $q_1+q_2+q_3+q_4=0$, $q_i^2=q_i$, 
$q_1^*+q_2^*+q_3^*+q_4^*=0$ and $(q_i^*)^2=q_i^*$ on a dense invariant domain
of ${\mathfrak A}$.
\end{remark}

\vspace{0.1cm}

{\bf 3. Again representations}

\vspace{0.1cm}

Next result 
 shows that the class of unbounded
representations of ${\mathcal A}_{4,0}$ satisfying the conditions of 
Definition~\ref{def} is $*$-wild (see \cite{t} for the definition of $*$-wild
unbounded representations). 
\begin{proposition}
The class of integrable representations of ${\mathcal A}_{4,0}$ is $*$-wild.
\end{proposition}
\begin{proof}
Let $\alpha, \beta>0$ and  let $\mathfrak{S}_2$ be the $*$-algebra generated 
by selfadjoint elements $a$ and $b$. Consider the set ${\mathfrak R}$ of all
representations $\pi$
of $\mathfrak{S}_2$   such that $||\pi(a)||\leq\alpha$, 
$||\pi(b)||\leq\beta$. Denote by
$\mathfrak{A}_{\alpha,\beta}$ the completion of
$\mathfrak{S}_2/\{z:|||z|||=0\}$ under $|||z|||=sup\{||\rho(z)||;\rho\in{\mathfrak R}\}$. 

Let $H$ be a separable infinite dimensional Hilbert space with an orthonormal 
basis $\{e_k\}_{k\in{\mathbb Z}}$, let $P_k$ be the orthoprojection onto
${\mathbb C}\langle e_k\rangle$, $k\in{\mathbb Z}$. We consider  operators $v$, $w$ defined by 
$ve_k=e_{k+1}$, 
$ve_{k+1}=e_k$ if $k$ is even and $we_k=e_{k+1}$, 
$we_{k+1}=e_k$ if $k$ is odd. Clearly, $(P_{2k}+P_{2k+1})H$ 
(respectively $(P_{2k+1}+P_{2k+2})H$) is invariant with respect to  $v$ (respectively $w$).

Let now 
\begin{eqnarray*}
&\tilde p=\sum_{k\ne 0}(-1)^{k+1}kP_k\otimes
\left(\begin{array}{cc}
e&0\\
0&e
\end{array}\right),\ 
\tilde q=\sum_{k\ne 0}(-1)^{k}kP_k\otimes
\left(\begin{array}{cc}
e&0\\
0&e
\end{array}\right),\\
&\tilde r= (\sum_{k\ne 0}(2k+1)vP_{2k}-2kvP_{2k+1})\otimes
\left(\begin{array}{cc}
e&0\\
0&e
\end{array}\right)+vP_0\otimes
\left(\begin{array}{ccc}
e&0&0\\
0&2e&0
\end{array}\right),\\
&\tilde s=(\sum_{k\ne 0}(2k+1)wP_{2k+2}-(2k+2)P_{2k+1})\otimes
\left(\begin{array}{cc}
e&0\\
0&e
\end{array}\right)+wP_0\otimes
\left(\begin{array}{ccc}
e&e&a+ib\\
0&e&e
\end{array}\right).
\end{eqnarray*}
Here $e$ is the identity element in $\mathcal{A}_{\alpha,\beta}$.
We write ${\mathcal H}$ for the Hilbert space 
$P_0H\oplus P_0H\oplus P_0H\oplus ((I-P_0)H\oplus (I-P_0)H)$ 
Let $CB({\mathcal H})$ be the $C^*$-algebra of compact operators on
${\mathcal H}$.
Direct verification shows that $\tilde p$, $\tilde q$, $\tilde r$, $\tilde s$ 
are affiliated with the $C^*$-algebra
$CB({\mathcal H})\otimes {\mathcal A}_{\alpha, \beta}$ (the completion of
the algebraic tensor product of $CB({\mathcal H})$ and 
${\mathcal A}_{\alpha, \beta}$ with respect to a $C^*$-norm, it does not 
depend which one). Moreover, since
any representation of $CB({\mathcal H})\otimes{\mathcal A}_{\alpha,\beta}$
is of the form $V^{-1}(id\otimes\pi)V$, where $V$ is a unitary operator,
$id$ is the identical representation of $CB({\mathcal H})$ and $\pi$ is
a representation of ${\mathcal A}_{\alpha,\beta}$,  one can show that
$\tilde p$, 
$\tilde q$, $\tilde r$, $\tilde s$, $\tilde r^*$, $\tilde s^*$   
 separate representations of 
$CB({\mathcal H})\otimes {\mathcal A}_{\alpha, \beta}$, i.e., if $\pi_1$, 
$\pi_2$ are different non-degenerate representations of $CB({\mathcal H})\otimes {\mathcal A}_{\alpha, \beta}$ then $\pi_1(x)\ne\pi_2(x)$, where $x$ is one of $p$, $q$, $r$, $s$.  
In fact, if $V^{-1}_1(id\otimes\pi_1)(x)V_1=V_2^{-1}(id\otimes\pi_2)(x)V_2$,
$x=p,q,r,s,r^*,s^*$, direct verification shows that $V_2V_1^{-1}=I\otimes V$,
where $V\pi_1=\pi_2V$ and therefore $V^{-1}_1(id\otimes\pi_1)V_1=
V_2^{-1}(id\otimes\pi_2)V_2$.
Besides,
since $(I+\tilde p^2)^{-1}=\sum_{k\ne 0}(1+k^2)^{-1}P_k\otimes  
\left(\begin{array}{cc}
e&0\\
0&e
\end{array}\right)$, $(I+\tilde p^2)^{-1}\in CB({\mathcal H})\otimes 
{\mathcal A}_{\alpha, \beta}$. Therefore, by \cite[Theorem~3.3]{wor2}, 
$\tilde p$,
 $\tilde q$, $\tilde r$, $\tilde s$
generate the $C^*$-algebra $CB({\mathcal H})\otimes {\mathcal A}_{\alpha, \beta}$.

Let $D=l.s.\ \{a\otimes b\mid a\in CB({\mathcal H}),
a\in{\mathcal F}, b\in {\mathcal A}_{\alpha,\beta}\}$, where ${\mathcal F}$ 
is the space of finite-dimensional operators in ${\mathcal H}$.
Then $D$ is dense in $CB({\mathcal H})\otimes {\mathcal A}_{\alpha, \beta}$
and invariant with respect to $\tilde p$, $\tilde q$, $\tilde r$, $\tilde s$,  
$D$ is a core
for the elements $\tilde p$, $\tilde q$, $\tilde r$, $\tilde s$ and 
$\tilde p$, $\tilde q$, $\tilde r$, $\tilde s$ satisfy relations 
(\ref{rel1})--(\ref{rel2})  
on $D$. 
Moreover, with $\psi(p)=\tilde p$, $\psi(q)=\tilde q$, $\psi(r)=\tilde r$, $\psi(s)=\tilde s$  the 
 representation ($\pi(\psi)(p))$, $\pi(\psi(q))$, $\pi(\psi(r)$,
$\pi(\psi(s))$) of ${\mathcal A}_{4,0}$ satisfies the condition of 
Definition~\ref{def}
for any  representation $\pi$ of
 $CB({\mathcal H})\otimes {\mathcal A}_{\alpha, \beta}$. From this it follows 
that
the class $R$ is $*$-wild.


The mapping $\psi$ defines a  functor $F_{\psi}$ from the category ${\mathcal Rep}({\mathcal A}_{\alpha,\beta})$
 of 
non-degenerated
representations 
of ${\mathcal A}_{\alpha,\beta}$ to the category
${\mathcal Rep}({\mathcal A}_{4,0})$ as follows:
\begin{itemize}
\item $F_{\psi}(\pi)(x)=(id\otimes\pi)(\psi(x))$ for any 
$\pi\in{\mathcal Rep}({\mathcal A}_{\alpha,\beta})$, $x=p,q,r,s$, 
\item $F_{\psi}(A)=E\otimes A$ for any operator $A$ intertwining $\pi_1$ and 
$\pi_2\in{\mathcal Rep}({\mathcal A}_{\alpha,\beta})$.
\end{itemize}
Since $id\otimes\pi$ is a representation of $CB({\mathcal H})\otimes {\mathcal A}_{\alpha, \beta}$ it can be uniquely extended to affiliated elements $\psi(p)$,
$\psi(q)$, $\psi(r)$, $\psi(s)$.
It follows from \cite{t} that
the functor  $F_{\psi}$ is full.
\end{proof}

\section{Representations of algebras $\Q_{n,\lambda}$  and $*$-algebras
$\Q_{n,\lambda}(*)$, $n\geq 5$}

\subsection{Algebras $\Q_{n,\lambda}$, $n\geq 5$, and their representations}
For each $n\geq 5$ and $\lambda\in {\mathbb C}$ the algebra $\Q_{n,\lambda}$
is non-zero and contains as a subalgebra the free algebra with two generators,
\cite{rss}.

In this paper we do not give the description of the whole set 
$\Lambda_{n,fd}$ for $n\geq 5$ but some facts concerning this set.
For other results see, for example, \cite{Wu1,wang}.

As it was noticed before, 
$\Lambda_{n,fd}\subset {\mathbb Q}$. On the other hand, $\Lambda_{n,fd}$
contains the set $\Sigma_{n,fd}=\{\alpha\in {\mathbb R}\mid
\exists H, dim H<\infty, \text{ orthoprojections } P_1,\ldots,P_n \text{ such
that } \sum P_k=\alpha I\}$, the last being studied in
\cite{krs}. By \cite{krs}, the following statement holds.
\begin{proposition}
$${\mathbb Q}\supset\Lambda_{n,fd}\supset\Lambda_n^1\cup\Lambda_n^2$$
where 
\begin{eqnarray*}
\Lambda_n^1=\{0, 1+\frac{1}{(n-1)},1+\frac{1}{(n-2)-
\displaystyle\frac{1}{(n-1)}},\ldots,1+\frac{1}{(n-2)-
\displaystyle\frac{1}{(n-2)-\frac{1}{\ddots-
\displaystyle\frac{1}{(n-1)}}}},\ldots\},\\
\Lambda_n^2=\{1,1+\frac{1}{(n-2)},1+\frac{1}{(n-2)-
\displaystyle\frac{1}{(n-2)}},\ldots,1+\frac{1}{(n-2)-
\displaystyle\frac{1}{(n-2)-\frac{1}{\ddots-
\displaystyle\frac{1}{(n-2)}}}},\ldots,\}
\end{eqnarray*}
\end{proposition} 
As to the description of finite-dimensional representations 
of $\Q_{n,\Lambda}$,
$\lambda\in\Lambda_{n,fd}$, up to similarity, this is an open question now.

Concerning the set $\Lambda_{n,bd}$, $n\geq 5$, we have the following
\begin{proposition}\label{cuntz}
$$\Lambda_{n,bd}={\mathbb C},\quad (n\leq 5).$$
\end{proposition}
\begin{proof}
For each $\lambda\in{\mathbb C}$, we give, following \cite{rs}, a concrete 
construction of five idempotents
$Q_i\in B(H)$, whose sum is equal to $\lambda I$.
Let $H=l_2\oplus l_2\oplus l_2$ and let $\I$ denote the identity 
operator on $l_2$. Define $Q_i$, $i=1,\ldots,5$, in the following way:
\begin{eqnarray*}
&Q_1=\left(\begin{array}{ccc}
a\I&3a\I&b\I\\
a\I&3a\I&b\I\\
a\I&3a\I&b\I
\end{array}\right), Q_2=\left(\begin{array}{ccc}
a\I&-3a\I&b\I\\
-a\I&3a\I&-b\I\\
a\I&-3a\I&b\I
\end{array}\right)\\
&Q_3=\left(\begin{array}{ccc}
4a\I&0&-2b\I\\
0&0&0\\
-2a\I&0&b\I
\end{array}\right),
Q_4=\left(\begin{array}{ccc}
2c\I&0&2dcS_1^*\\
0&2c\I&2dcS_2^*\\
S_1&S_2\I&d\I
\end{array}\right)\\
&Q_5=\left(\begin{array}{ccc}
2c\I&0&-2dcS_1^*\\
0&2c\I&-2dcS_2^*\\
-S_1&-S_2\I&d\I
\end{array}\right)
\end{eqnarray*} 
where $a=(5-2\lambda)/6$, $b=(4\lambda-7)/3$, $c=(3\lambda-5)/4$, 
$d=(7-3\lambda)/2$, and $S_1$, $S_2$ are operators of  a representation of 
the Cuntz algebra
${\mathcal O}_2$ (\cite{cuntz}). Recall that ${\mathcal O}_2$ is a unital 
$*$-algebra generated by
$s_1$, $s_2$, $s_1^*$, $s_2^*$ and relations $s_1^*s_2=0$, 
$s_1^*s_1=e=s_2^*s_2=s_1s_1^*+s_2s_2^*$.
Direct verification shows that $Q_i=Q_i^2$, $i=1,\ldots,5$, and 
 $Q_1+Q_2+Q_3+Q_4+Q_5=\lambda I$.
\end{proof}
Note that the  construction which is given in the proposition
is a generalization of an example in \cite{bes} of
five idempotents with zero sum.

\subsection{$*$-Algebras $\Q_{n,\lambda}(*)$, $n\geq 5$ and their
$*$-representations by bound\-ed operators}
As we already know, $\Lambda_{n,bd}={\mathbb C}$ if $n\geq 5$.
As to representation of $\Q_{n,\lambda}(*)$ we have the following 
\begin{proposition}
$*$-Algebra $\Q_{n,\lambda}(*)$ is not of type $I$ for each 
$\lambda\in{\mathbb C}$ and $n\geq 5$,
 i.e.\ for each $\lambda\in{\mathbb C}$ it has a factor-representation which 
is not of type $I$.
\end{proposition}
\begin{proof}
It is enough to show that the $*$-algebras $\Q_{5,\lambda}(*)$ is not of type
$I$ for each $\lambda\in{\mathbb C}$. Assume first that $\lambda\ne 2$.
It is known that
the Cuntz algebra ${\mathcal O}_2$ is of not type $I$. So, there exists
a factor-representation, $\rho$, of ${\mathcal O}_2$ such that the double
commutant $\rho({\mathcal O}_2)''$ is not of type $I$.
Consider now a representation, $\pi$,  of $\Q_{5,\lambda}(*)$ given in 
Proposition~\ref{cuntz} with $S_i=\rho(s_i)$, $i=1,\ldots,5$.
Direct calculations show that the commutant $\pi(\Q_{5,\lambda}(*))'$ coincides
with $\{diag(C,C,C)\mid C\in\rho({\mathcal O}_2)'\}$
and $\pi(\Q_{5,\lambda}(*))''=M_3({\mathcal N})$, where 
${\mathcal N}=\rho({\mathcal O}_2)''$. Since ${\mathcal N}$ is not of type I,
$M_3({\mathcal N})$ is not of type $I$, completing the proof.

If $\lambda=2$, then $\Q_{5,\lambda}(*)$ is $*$-wild by 
Proposition~\ref{wild1} and
therefore is not of type $I$.
\end{proof}
For many $\alpha\in{\mathbb R}$ we can say even more:
there exist $\alpha\in{\mathbb R}$ such that $\Q_{n,\lambda}(*)$ and even
${\mathcal P}_{n,\alpha}$ is $*$-wild. 

\begin{proposition}
The $*$-algebras $\Q_{n,\alpha}(*)$ ($n\geq 5$) are $*$-wild for
$\alpha$ from the following sets:

(a) $\Lambda_{4,bd}=\{2\pm 2/k (k\in{\mathbb N}), 2\}$,

(b) $\Lambda_{n,orb(2)}=\{\alpha_0=2,\alpha_k=(n-1)-1/(\alpha_{k-1}-1), 
k\in{\mathbb Z}\}$,

(c) $\Lambda_{n,orb(n/2)}=\{\alpha_0=n/2,\alpha_k=(n-1)-1/(\alpha_{k-1}-1), 
k\in{\mathbb Z}\}$.
\end{proposition}
\begin{proof}
(a) $\Q_{n,\alpha}(*)$ ($n\geq 5$) is $*$-wild for 
$\alpha\in\Lambda_{4,bd}$ because $Q_{4,\alpha}(*)$ is $*$-wild for every
$\alpha\in\Lambda_{4,bd}$ by Proposition~\ref{wild2}.

(b) By \cite{os}[Theorem~57], \cite{krusam}[Theorem~4] the unital $*$-algebra
${\mathcal P}_{3,\perp 2}={\mathbb C}\langle r,r_1,r_2\mid
r=r^*,r^2=r,r_i^*=r_i,r_i^2=r_i, r_1r_2=0\rangle$ is $*$-wild.
Setting $\psi(p_1)=r$, $\psi(p_2)=e-r$, $\psi(p_3)=r_1$, $\psi(p_4)=r_2$,
$\psi(p_5)=e-r_1-r_2$ for the generators $p_i$, $i=1,\ldots,5$, of 
${\mathcal P}_{5,2}$, we obtain a $*$-epimorphism from ${\mathcal P}_{5,2}$
to ${\mathcal P}_{3,\perp 2}$ so that  ${\mathcal P}_{3,\perp 2}$ is
a factor $*$-algebra of ${\mathcal P}_{5,2}$. This shows that 
${\mathcal P}_{5,2}$ and therefore ${\mathcal P}_{n,2}$, $n\geq 5$, is 
$*$-wild.
Then the  Coxter functors $F:{\mathcal P}_{n,\alpha}\to
{\mathcal P}_{n,1+1/(n-\alpha-1)}$ and $R:{\mathcal P}_{n,\alpha}\to
{\mathcal P}_{n,n-1-1/(\alpha-1)}$, constructed in \cite{krs},
spread out the $*$-wildness to all points $$\alpha\in\Lambda_{n,orb(2)}=\{
\ldots,\alpha_{-1}=1+\frac{1}{n-3},\alpha_0=2,\alpha_1=n-2,
\ldots,\alpha_k=(n-1)-\frac{1}{\alpha_{k-1}-1},\ldots\}.$$
The set $\Lambda_{n,orb(2)}$ is the two-sided orbit of the point $\{2\}$
with respect to the dynamical system $\alpha\to f(\alpha)=(n-1)-1/(\alpha-1)$.

(c) By \cite{os}[Theorem~55],  the unital
$*$-algebra  
${\mathcal P}_{3,2anti}={\mathbb C}\langle w_i, i=1,2,3\mid
w_i^*=w_i,w_i^2=e,w_1w_2=w_2w_1=0\rangle$ is $*$-wild. 
The $*$-algebra
${\mathcal P}_{5,5/2}$ is itsfactor $*$-algebra with
corresponding $*$-epimorphism $\psi$ given by
$\psi(p_1)=(w_1+\sqrt{3}w_2+2e)/4$, $\psi(p_2)=(w_1-\sqrt{3}w_2+2e)/4$,
$\psi(p_3)=(-w_1+e)/2$,
$\psi(p_4)=(w_3+e)/2$, $\psi(p_5)=(-w_3+e)/2$.
Since ${\mathcal P}_{5,5/2}$ is a factor-$*$-algebra
of ${\mathcal P}_{2n+1,(2n+1)/2}$ for each $n\geq 5$, we obtain that
${\mathcal P}_{2n+1,(2n+1)/2}$ is $*$-wild for any $n\geq2$. 
Since the $*$-algebra ${\mathcal P}_{5,2}$ is $*$-wild, by the same arguments,
we obtain that ${\mathcal P}_{6,3}$  and ${\mathcal P}_{2n,n}$, $n\geq 3$
are $*$-wild. Then, using the Coxter functors $F$ and $R$ we get that
the $*$-algebras ${\mathcal P}_{n,\alpha}$ are $*$-wild for any
$$\alpha\in\Lambda_{n,orb(n/2)}=\{\ldots,\alpha_{-1}=1+\frac{2}{n-2},
\alpha_0=\frac{n}{2}, \ldots, \alpha_k=(n-1)-
\frac{1}{\alpha_{k-1}-1},\ldots\}.$$ 
The set $\Lambda_{n,orb(n/2)}$ is the two-sided orbit of the point $\{n/2\}$
with respect to the dynamical system $\alpha\to f(\alpha)=(n-1)-1/(\alpha-1)$.
\end{proof}

Restricting ourselves to $*$-representations of ${\mathcal P}_{n,\alpha}$,
or, equivalently, $*$-representations of $\Q_{n,\alpha}(*)$ with the condition
that the images of $q_i$ are selfadjoint, we can give the  full classification
of such $*$-representations for $\alpha\in\Lambda_n^1\cup\Lambda_n^2$.
If $\alpha\in\Lambda_n^1$ there exists a unique, up to unitary equivalence, irreducible 
representation of ${\mathcal P}_{n,\alpha}$, but if $\alpha\in\Lambda_n^2$ 
there 
are $n$ unitarily non-equivalent irreducible representations of
${\mathcal P}_{n,\alpha}$ (see \cite{krs}).   
\subsection{Representations of $\Q_{n,\lambda}$, $n\geq 5$, by unbounded
operators}
We do not study here unbounded representations of $\Q_{n,\lambda}$, 
$n\geq 5$ as the structure of bounded representations of 
$\Q_{n,\lambda}$ ($n\geq 5$) is already very complicated.

\vspace{1cm}

\noindent
{\bf Acknowledgements.}
The work was completed during the visit of the first  author
Chalmers University of Technology and G\"oteborg University
 in May-June 2001. The author is pleased to thank the department there for 
their warm hospitality. This work was supported by the Swedish Royal Academy 
of Sciences.


\begin{thebibliography}{50}
\bibitem[BES]{bes} H.Bart, T.Ehrhardt and B.Silbermann, {\it Zero sums of
idempotents in Banach algebras}. Integr. Equat. Oper. Th. {\bf 19} (1994),
125-134.

\bibitem[Ba]{bax}
{\it R.J.~Baxter},  Exactly solved models in statistical mechanics. 
Academic Press, Inc. [Harcourt Brace Jovanovich,
Publishers], London.

\bibitem[Bo]{bondar}
V.M.~Bondarenko,  {\it On certain wild algebras generated by idempotents}.
 Methods Funct. Anal. Topology {\bf 5} (1999), no. 3,
1--3. 
\bibitem[BGKKRSS]{bgkkrss}
A.~B\"ottcher,  I.~Gohberg, Yu.~Karlovich,  N.~Krupnik, S.~Roch, 
B.~Silbermann, I.~Spitkovsky, {\it Banach algebras generated by
$N$ idempotents and applications}.  Oper. Theory Adv. Appl., 90,
Birkhäuser, Basel, (1996), 19-54.
\bibitem[Cu]{cuntz}
J.Cuntz, {\it Simple $C^*$-algebras generated by isometries}. Comm. Math. Phys.
{\bf 57}, (1977), 173-185.
\bibitem[D]{yugoslav}
D.\u {Z}.~Dokovi\'c, {\it Unitary similarity of projectors}.
 Aequationes Math. {\bf 42} (1991), no. 2-3, 220--224.
\bibitem[DF]{df}
{\it P.~Donovan, M.R.~Freislich},  The representation theory of finite 
graphs and associated algebras. Carleton
Mathematical Lecture Notes, No. 5. Carleton University, Ottawa, Ont., 1973. 
\bibitem[EK]{ek} 
{\it D.E.~Evans,  Yu.~Kawahigashi}, Quantum symmetries on operator algebras. 
Oxford Mathematical Monographs.
Oxford Science Publications. The Clarendon Press, Oxford University Press, 
New York, 1998.
\bibitem[ERSS]{erss}
T.~Ehrhardt, V.~Rabanovich, Yu.~Samo\u\i{}lenko, B.~Silbermann,
{\it On decomposition of the identity into a sum of idempotents}. Methods
 Funct. Anal. Topology, {\bf 7} (2001), no. 2.
\bibitem[GGV]{ggv}
{\em I.M.~Gel'fand, M.I.~Graev and N.Ya.~Vilenkin}, Generalized functions.
vol. 5, Academic Press, New York and London, 1966.
\bibitem[GP]{gp}
I.M.~Gel'fand, V.A.~Ponomarev, {\it Quadruples of subspaces of a 
finite-dimensional space}. Dokl.\ Acad.\ Nauk SSSR, {\bf 197} (1971), no.4,
762-765, (Russian).
\bibitem[GoP]{gorp} M.F.~Gorodni\u\i, G.B.~Podkolzin, {\it Irreducible 
representations of graded Lie algebras}. (Russian) Spectral theory of 
operators
and infinite-dimensional analysis, 66--77, iii, Akad. Nauk Ukrain. SSR, Inst. Mat., Kiev, 1984. 
\bibitem[In]{inoue}
{\it A.~Inoue}, Tomita-Takesaki theory in algebras of unbounded operators.
 Lecture Notes in Mathematics, 1699.
Springer-Verlag, Berlin, 1998.
\bibitem[KS2]{krusam1}
S.~Krugliak,  Yu.~Samo\u\i{}lenko, {\it Structure theorems for families of idempotents}, Ukra\"in. Mat. Zh.
{\bf 50} (1998), no.~4, 523-533, (Russian).  
\bibitem[KS2]{krusam}
S.~Krugliak,  Yu.~Samo\u\i{}lenko, {\it On complexity of description of 
representations of $*$-algebras generated by
idempotents}. Proc. Amer. Math. Soc. {\bf 128} (2000), no. 6, 1655--1664. 
\bibitem[KRS]{krs}
S.~Krugliak, V.~Rabanovich, Yu.~Samo\u\i{}lenko, {\it On sums of projections}.
(in preparation).
\bibitem[Na]{naz}
L.A.~Nazarova,  {\it Representations of a tetrad}. Izv. Akad. 
Nauk SSSR Ser. Mat. {\bf 31} (1967), 1361--1378, (Russian).
\bibitem[N]{nel}
E. Nelson, {\it Analytic vectors}. Ann. Math., {\bf70}, no.~3, 572--615 (1959).
\bibitem[OS1]{osb} 
V.~L. Ostrovsky\u\i{}, Yu.~S. Samo\u\i{}lenko, {\it Unbounded operators
satisfying non-{L}ie commutation relations},
 Repts. Math. Phys. {\bf 28}, 1, (1989), 91--103.
\bibitem[OS2]{os}
 {\em Ostrovski\u\i{} V., 
   Samo\u\i{}lenko Yu.}, Introduction to the Theory 
   Representation of Finitely Presented 
   $*-$algebras.
   1. Representations by bounded operators.
   The Gordon and Breach Publ. Group, London, 1999.
\bibitem[P]{popovich}
S.~Popovych, {\it Unbounded idempotents}. Methods Funct. Anal. Topology 
{\bf 5} (1999), no. 1, 95--103. 
\bibitem[RS]{rs}
V.~Rabanovich, Yu.~Samo\u\i{}lenko, {\it When sum of idempotents or 
projections is a multiple of unity}.   
Funktsional. Anal. i Prilozhen. {\bf 34} (2000), no. 4, 91--93.
\bibitem[RSS]{rss}
V.I.~Rabanovich, Yu.S.~Samo\u\i{}lenko, A.V.~Strelets,{\it On identities in
algebras $\Q_{n,\lambda}$ generated by idempotents}. Ukra\"in. Mat. Zh.
{\bf 53} (2001), no.8.
\bibitem[RaS]{ras}
 I.~Raeburn, A.M.~Sinclair,  {\it The $C\sp *$-algebra generated by two 
projections}. Math. Scand. {\bf 65} (1989), no. 2,
278--290.
\bibitem[ST]{samtur}
Yu.S.~Samo\u\i lenko, L.B. Turovskaya, {\it On representations of 
$*$-algebras by unbounded operators}, (Russian)
Funkt. Anal. i Prilozhen. 31 (1997), no. 4,80--83 translation in 
Funct. Anal. Appl. 31 (1997), no. 4, 289--291 (1998).
\bibitem[S]{Shbook} 
{\it K.~Schm\"udgen}, Unbounded operator algebras and representation theory.
Birkh{\"a}user Verlag, Berlin, 1990.
\bibitem[T]{t} L.~Turowska, {On the complexity of $*$-algebras representations 
by unbounded operators}. Preprint, Department of Mathematics, Chalmers 
University of Technology and G\"oteborg University, 2000, pp.21, to appear
in Proc.Amer.Math.Soc.
\bibitem[V]{V} {\em N.Ya.~Vilenkin},
Special functions and the theory of group representations. Second edition.
 ``Nauka'', Moscow, 1991 (Russian). 
\bibitem[Wa]{wang}
Jin-Hsien Wang, {\it The length problem for a sum of idempotents}.
 Linear Algebra Appl. {\bf 215} (1995), 135--159.
\bibitem[W1]{wor1} S.L. Woronowicz, {\it Unbounded elements affilated
with $C^*$-algebras and no-compact quantum groups}. Commun.Math.Phys. 136
(1991), 399-432.
\bibitem[W2]{wor2} S.~L. Woronowicz, \emph{{$C^*$}-algebras generated by 
unbounded elements}.
  Reviews in Mathematical Physics \textbf{7} (1995), no.~3, 481--521.
\bibitem[WN]{wn}
S.~L. Woronowicz and K. Napiorkowski,\emph{Operator theory in the $C^*$-algebra
framework}. Rep. Math. Phys. \textbf{31} (1992), 353-371.  
\bibitem[Wu1]{Wu}
P.Y.~Wu, {\it Additive combination of special operators}. Funct.\ Anal.\ and
Oper.\ Theory, Banach Center Publ., Warszawa, {\bf 30} (1994), 337-361.
\bibitem[Wu2]{Wu1}
P.Y.~Wu,  {\it Sums of idempotent matrices}. Linear Algebra Appl. 142 (1990), 
43--54. 
\end{thebibliography}
\end{document}